%% file: main.tex
\documentclass[a4paper,final]{article}
\usepackage{amsmath}
\usepackage{amsthm}
\usepackage{amssymb}
\usepackage{graphicx}
\usepackage{multirow}
\usepackage{longtable}
\usepackage{wrapfig}
\usepackage{float}
\usepackage{xcolor}
\usepackage{hyperref}
\usepackage{setspace}
\usepackage{mathtools}
\usepackage{algorithm}
\usepackage{algorithmic}
\renewcommand{\algorithmiccomment}[1]{\bgroup\hfill//~#1\egroup}
\usepackage{eurosym}
\usepackage[utf8]{inputenc}
\usepackage[english]{babel}
\usepackage{amsfonts}
\usepackage{authblk}
\usepackage{xcolor}
\usepackage[normalem]{ulem}
\newtheorem{theorem}{Theorem}

\newtheorem{definition}[theorem]{Definition}

\newtheorem{lemma}[theorem]{Lemma}

\newtheorem{proposition}[theorem]{Proposition}

\title{One Benders cut to rule all schedules in the neighbourhood}

\author[1,*]{Ioannis Avgerinos}
\author[1]{Ioannis Mourtos}
\author[1]{Stavros Vatikiotis}
\author[1]{Georgios Zois}
\affil[1]{\small ELTRUN Lab, Department of Management Science and Technology, Athens University of Economics and Business,\protect\\ \small Patision 76, 10434, Athens, Greece E-mail: {\tt iavgerinos@aueb.gr, mourtos@aueb.gr, stvatikiotis@aueb.gr, georzois@aueb.gr}}

\affil[*]{Corresponding Author} 
\date{} 

\begin{document}
\maketitle

\begin{abstract}
Logic-Based Benders Decomposition (LBBD) and its Branch-and-Cut variant, namely Branch-and-Check, enjoy an extensive applicability on a broad variety of problems, including scheduling. Although LBBD offers problem-specific cuts to impose tighter dual bounds, its application to resource-constrained scheduling remains less explored. Given a position-based Mixed-Integer Linear Programming (MILP) formulation for scheduling on unrelated parallel machines, we notice that certain $k-$OPT neighbourhoods could implicitly be explored by regular local search operators, thus allowing us to integrate Local Branching into Branch-and-Check schemes. After enumerating such neighbourhoods and obtaining their local optima - hence, proving that they are suboptimal - a local branching cut (applied as a Benders cut) eliminates all their solutions at once, thus avoiding an overload of the master problem with thousands of Benders cuts. However, to guarantee convergence to optimality, the constructed neighbourhood should be exhaustively explored, hence this time-consuming procedure must be accelerated by domination rules or selectively implemented on nodes which are more likely to reduce the optimality gap. In this study, the realisation of this idea is limited on the common 'internal (job) swaps' to construct formulation-specific $4$-OPT neighbourhoods. Nonetheless, the experimentation on two challenging scheduling problems (i.e., the minimisation of total completion times and the minimisation of total tardiness on unrelated machines with sequence-dependent and resource-constrained setups) shows that the proposed methodology offers considerable reductions of optimality gaps or faster convergence to optimality. The simplicity of our approach allows its transferability to other neighbourhoods and different sequencing optimisation problems, hence providing a promising prospect to improve Branch-and-Check methods.
\end{abstract}


\input{introduction.tex}
\input{framework.tex}
\input{scheduling.tex}
\input{local_branching.tex}
\input{experiments.tex}

\hfill \break
\textbf{Acknowledgments. }This work was supported by MODAPTO Horizon 2020 [grant number 101091996] (\url{https://modapto.eu/}.

\setstretch{1.05}
\input{references.tex}

\newpage
\end{document}

%% file: introduction.tex
\section{Introduction}
    \subsection{Context}
The industrial importance of production scheduling has been motivating the application of optimisation methods on production settings of increasing diversity. As this class of optimisation problems is known to be $\mathcal{NP}$-hard \cite{Len77}, exact methods typically require a dramatically large run-time to reach optimality even for instances of moderate scale. This explains why the largest part of related literature examines (meta-)heuristics that exchange the quality of the obtained solution with a faster computation.

Although machine scheduling problems appear daunting for an MILP solver, decomposing such problems raises some hope for an exact method. \emph{Benders Decomposition} \cite{benders} exploits the structure of an optimisation problem to partition it into the \emph{master} problem that provides dual bounds and the \emph{subproblem} that provides primal bounds. The two exchange knowledge in the form of \emph{Benders cuts} that are iteratively added to the master problem. For the classical variant of Benders Decomposition, these cuts are derived from duality theory on real variables. The \emph{Logic-Based} (LBBD) variant \cite{hooker03} allows for cuts also on integer variables. As total run-time of LBBD algorithms is spent mostly on the master problem (e.g. see \cite{beck}) the hybridisation of LBBD with \emph{Branch-and-Cut}, known as \emph{Branch-and-Check} \cite{thor}, is efficient for partitions having a quickly solvable subproblem \cite{beck}. This is our starting point.
    
In this paper, we consider a partitioning scheme for scheduling on unrelated parallel machines with sequence-dependent setup times and a renewable resource restricting the number of simultaneous setups. As in Branch-and-Check, each new integer solution found for the master problem activates the solution of a subproblem, thus obtaining a new Benders cut that eliminates suboptimal solutions. Our idea is to improve performance by generating, per integer solution of the master problem, a local branching constraint that eliminates all solutions in a standard neighbourhood surrounding that solution of the subproblem. Contrary to typical local branching for LBBD (e.g., \cite{jeihoonian}), this single cut is produced through local search that produces all solutions in a given neighbourhood and then solves a single Constraint Programming (CP) model to select the best among them. Our approach further strengthens this modified LBBD with simple domination rules to speed up the exploration of neighbourhoods and a class of strong valid inequalities tightening the lower bounds derived by the master problem. This allows us to tackle a large class of scheduling problems under two, rather challenging, objectives: total completion time and total tardiness. 

\subsection{Relevant literature}
        The concept of generating strong(er) cuts to accelerate convergence to optimality is as old as the \emph{Cutting-Planes} algorithm itself: \cite{gomory} states that the globally valid, yet usually ineffective, cuts could be replaced with ``\emph{particularly effective}" ones, exploiting the case-specific structure of the optimisation problem to be solved. \cite{benders} admits that ``\emph{adding more than one or two"} cuts at each iteration may imply a reduction of the total number steps which are required for convergence to be reached, always considering that this should occur in an effective way, ``\emph{avoiding the addition of redundant constraints}". \cite{magnati} relate the selection of proper formulations and cuts with faster performance of Benders Decomposition algorithms.
                
Recently, \cite{asl} showed that obtaining multiple solutions of the subproblem per iteration, hence generating multiple \emph{Benders cuts} per iteration, reduces the number of solutions of the master problem, in which 90\% of total solution time is spent on average. Also, \cite{seo} showed a scheme for the selection of the `closest' optimality cut, to accelerate the performance of Benders decomposition. 
\cite{jain} present a set of globally valid cuts and recommend using problem-specific strong cuts to enhance LBBD performance. We are particularly motivated by the ``\emph{strong cuts}" of \cite{harjunkoski}, i.e., inequalities eliminating sets of infeasible or suboptimal solutions, rather than eliminating a single solution at each step. 
        
The concept of hybrid methods, which combine the accuracy of exact algorithms and the speed of (meta-)heuristics, is also relevant (see the review in \cite{Puchinger05}). Solving the master problem or the subproblem of an LBBD heuristically implies faster computation but compromises convergence to optimality. As highlighted by \cite{hooker19}, obtaining suboptimal upper bounds (for minimisation problems) from the subproblem abolishes the validity of LBBD, as the generated cuts may eliminate the optimal solution. Regarding however the master problem, each solution defines a lower bound of the minimum objective value, hence a premature termination of the master problem retains validity, as long as the best dual bound is obtained at the end \cite{tran}. Indicatively, \cite{raidl14, raidl15} integrate meta-heuristics within LBBD to provide better solutions but without guaranteeing optimality.

The concept of Local Branching has already been used within Benders Decomposition by \cite{rei09}. Loosely described, given a solution of the master problem, the subproblem is solved iteratively for all, or at least for a sufficiently large number, of the solutions of a local neighbourhood. Each newly-obtained upper bound triggers the addition of a new Benders cut to the master problem, hence at each iteration, multiple suboptimal (or infeasible) solutions are eliminated. The large number of cuts, which are inevitably derived for problems of practical scale, requires a cut-management strategy to ensure the solvability of the master problem. Duality properties could be exploited to obtain the strongest cut, although this approach is mostly relevant to the classical Benders with subproblems of only continuous variables, as in \cite{zhu} for a scheduling problem and in \cite{jeihoonian} for a network design problem.
        
Our work aims not only at improving LBBD through stronger cuts arising from local search but also on providing an approach for all known objectives of parallel machine scheduling (makespan, sum of completion times, tardiness) and under diverse settings, e.g., uniform or unrelated machines, setup times that depend or not on the sequence, unlimited or limited resources. Hence, let us take a quick look also on this part of the literature. Although approximation algorithms have received strong interest \cite{Hall97}, exact methods relying upon LBBD are becoming more popular. The majority of current methods concerns the minimisation of the makespan, as it apparently allows for strong lower bounds \cite{Pey19, Pey20}. We have also combined heuristics and LBBD for makespan minimisation of splittable jobs on unrelated machines \cite{ijpr2022}. The work of \cite{tran} presents a Branch-and-Check method that partitions the problem into the \emph{assignment phase} of jobs to machines (master problem), and the \emph{sequencing phase} of jobs into the same machine (subproblem). 

Exact methods still struggle to solve the minimisation of total completion times and due-date related objectives. Although valid MILP formulations are known (e.g., \cite{Kim20, Mae20}), their evaluation typically concerns instances of small scale. This occurs due to the incapability of regular linear constraints to express the sequencing phase of jobs in an effective way, i.e., big-M constraints are known to provide weak lower bounds. This is partially addressed by the dominance and other properties of optimal schedules, as in \cite{Yal00} and \cite{Liaw03}. Regarding the minimisation of total completion time, the Benders Decomposition of \cite{Bul17} is a major contribution, as instances of considerable scale are tackled, without however resource constraints and sequence-dependent setup times. To the contrary, the minimisation of due-date related objectives implies additional challenges. Given a completion time $C_{j}$ for a job $j$ of due-date $d_{j}$, the respective value of tardiness is obtained by the nonlinear expression $T_{j} = \max\{0, C_{j} - d_{j}\}$. The linearisation of this expression is trivial but occurs in a large computational cost. The most significant progress on exact methods is observed by the Column Generation scheme of \cite{Pes22}, concerning the minimisation of weighted tardiness but on a single machine.

\subsection{Contribution}
        We combine local search with LBBD by introducing local branching cuts as Benders cuts that rule out all solutions in a neighbourhood. We show the virtue of this approach to  Parallel Machines Scheduling Problem, using a new decomposition that consists of a position-based formulation for the master problem \cite{Bektas} and results in a simpler subproblem. We strengthen that subproblem with a local branching cut, derived by considering regular local search operators for the solution of the master problem and prove that an associated neighbourhood is \textit{k}-OPT in local-branching terms \cite{fischetti03}. 
        
    Our idea extends a number of LBBD-accelerating concepts. The implementation of \emph{Local Branching} on Benders decomposition shows up in \cite{rei09} but restricted in the classical variant where the subproblem has only continuous variables. We extend this for LBBD, i.e., for mixed-integer subproblems. Most importantly, we no longer bother with the selection of the closest Benders cut, typically derived from duality theory \cite{seo}, as we opt for the generation of a single `supercut', eliminating an entire $k$-OPT neighbourhood per iteration. 
        
    Our LBBD partitioning scheme also has some novelties, e.g., once compared to \cite{tran} we notice that our approach has resource constraints hence not exploiting the improved cuts in \cite{tran}. Technically, a typical LBBD (e.g., \cite{tran}) solves the assignment of jobs to machines in the master problem and the sequencing of jobs assigned to machines in the subproblem. Our master problem does both but neglects resource-availability, leaving it to the subproblem. This implies a harder master problem and an easier subproblem, which is however enlarged and strengthened through local branching, and is more efficient in a Branch-and-Check context. We also examine ways to accelerate the exploration of the neighborhood and build all these as an LBBD code for minimisation of regular scheduling objectives. Further, our work focuses not on the well-studied makespan but on total tardiness and total completion times. Regarding the latter objective, it moves beyond \cite{Bul17} by tackling sequence-dependent and resource-constrained setups.

To evaluate our approach, we experiment on randomly generated instances of diverse features, spanning to a scale of $100$ jobs on $10$ machines. The results confirm that the proposed local branching scheme outperforms regular LBBD algorithms in terms of optimality gaps. In that regard, we consider our contribution substantial for the problem itself, as exact methods fail in even simple scheduling variants under total tardiness \cite{Kim20, Mae20} or cannot step up easily beyond a single machine \cite{Pes22}.
        
The remainder of this paper is structured as follows. Section \ref{sec:theory} recalls basics on LBBD and local branching to then present our variant. Our decomposition for parallel machine scheduling is presented in Section \ref{section:problem}. Section \ref{sec:lb} strengthens that decomposition, also by applying the ideas of Section \ref{sec:theory}. Our computational evidence regarding the minimisation of total completion times and total tardiness is discussed in Section \ref{sec:exp}. 

%% file: framework.tex
\section{Theory} \label{sec:theory}
\subsection{Logic-Based Benders decomposition}
        Consider an optimisation problem $\mathcal{P}$ defined on two vectors of variables $x,y$ as 
        
        $$\mathcal{P}: min\{f(x)+g(y)|C_1(x), C_2(y), C_3(x,y), x\in D_x, y\in D_y\},$$ 
        where functions $f,g$ and constraint sets $C_1,C_2,C_3$ are of arbitrary form. LBBD partitions $\mathcal{P}$ into a relaxation called the \emph{master problem} $\mathcal{M}$, which involves only variables $x,$ function $f(x)$ and constraints $C_1(x)$, leaving all other aspects of $\mathcal{P}$ to the \emph{subproblem} $\mathcal{S}$. Formally, $$\mathcal{M}: min\{z\geq f(x)|C_1(x), x\in D_x\},$$ and, given an optimal solution $\bar{x}$ of $\mathcal{M}$,  $$\mathcal{S}: min\{\zeta\geq f(\bar{x}) + g(y)|C_2(y), C_3(\bar{x}, y), y\in D_y\}.$$
        
        \noindent Therefore, if $\bar{z}$ is the value of an optimal solution of $\mathcal{P}$, $z\leq \bar{z}\leq \zeta$, because a solution to $\mathcal{S}$ is feasible for the original problem $\mathcal{P}$ and $\mathcal{M}$ is a relaxation of $\mathcal{P}$. As this optimality gap is typically large, LBBD proceeds in an iterative fashion by adding optimality cuts (also called Benders cuts) to $\mathcal{M}$ at each iteration $t$, with the obvious aim to increase the lower bound $z^{t}$ and decrease the upper bound $\zeta^{t}$. 
        
        Let $\bar{x}^{t}$ be the optimal solution to $\mathcal{M}$ at iteration $t$. To ensure convergence to an optimal solution of $\mathcal{M}$, LBBD requires a bounding function $B_{t}(\textsc{x})$ with two properties.
        \begin{description}
            \item{$P.1.$ }$B_{t}(x)\leq f(x) + g(y)$ for all feasible solutions $x,y$, and
            \item{$P.2.$ }$B_{t}(x) = \zeta^{t}$ if $x = \bar{x}^{t}$.
        \end{description}
        \noindent The typical bounding function (see for example \cite{hooker07}) is 
        \begin{equation}
            B_{t}(x) = \begin{cases}
                                    \zeta^{t} \quad \text{if }x = \bar{x}^{t}, \\
                                    0 \quad \text{otherwise},
                        \end{cases}
            \label{eq:Bounding_Function}
        \end{equation}  
        \noindent and triggers the generation of an optimality cut per iteration $t$, namely:
        \begin{equation}
            z\geq B_{t}(x). \label{eq:benders_cut}
        \end{equation}
        
        \noindent Adding (\ref{eq:benders_cut}) to the master problem at iteration $t$ yields that $z\geq \zeta^{t}$ if $x = \bar{x}^{t}$; otherwise, $z\geq f(x)$, as dictated by the original formulation of $\mathcal{M}$, yielding that (\ref{eq:Bounding_Function}) complies with properties $P.1$ and $P.2$. Also, $z^{t}\leq \bar{z}\leq \zeta^{t}$ at each iteration $t$ of the LBBD algorithm. By \cite[Theorem 1]{hooker07}, if the bounding function adheres to $P.1$ and $P.2$, and if the domain $D_{y}$ if finite, then the LBBD algorithm converges to the optimal solution of $\mathcal{P}$ after finitely many steps.

Solving $\mathcal{M}$ to optimality (or even to near optimality) at each iteration means that the convergence is slow, particularly if $\mathcal{M}$ is quite heavy. To the contrary, the Branch-and-Check approach solves the master problem only once. Whenever a new integer solution to $\mathcal{M}$  is found, the respective subproblem checks whether the obtained solution yields an improved primal bound. Most importantly, the subproblem also generates a Benders cut to be added to the continuing solution procedure of $\mathcal{M}$. Therefore, the index $t$ denotes in Branch-and-Check the counter of explored integer solutions and no longer the iteration. The convergence to optimality is ensured, as long as Properties $P.1$ and $P.2$ are valid for the generated cut. 
    
    \subsection{Local branching as Benders cuts} \label{section:local_branching}

Let $\bar{x}^{t},$ i.e., the $t^{th}$ integer solution of $\mathcal{M}$ found, be the \emph{reference solution} in terms of \cite{fischetti03}. Such a solution defines a $k$-OPT neighbourhood $N(\bar{x}^{t}, k)$ that includes all solutions $x$ having a symmetric difference at most $k$ with $\bar{x}^t$; i.e., at most $k$ variables receive in any $x \in N(\bar{x}^{t}, k)$ different values than in $\bar{x}^{t}$. 
        
Our twisted LBBD obtains a subproblem $\mathcal{S'}$ not by just fixing $x$ to the solution $\bar{x}^t$ of $\mathcal{M}$ but by simply adding to $\mathcal{P}$ the standard local branching constraint
                \begin{equation}
                 \sum_{j: \bar{x}_j^{t}=1}(1 - x_{j}) + \sum_{j: \bar{x}_j^{t}=0}x_{j}\leq k \label{eq:local_branch_cut}.
                \end{equation}
        \noindent Formally, $$\mathcal{S}'= min\{\eta\geq f(x)+g(y)|C_1(x), C_2(y), C_3(x,y), (\ref{eq:local_branch_cut}),x\in D_x,y\in D_y\}.$$ 
        
        \noindent That is, $\mathcal{S}'$ provides the optimal among all solutions of the $k$-OPT neighbourhood $N(\bar{x}^{t}, k)$, thus computing a bounding function $B'_{t}(x) \leq B_{t}(x)$, defined as
        \begin{equation}
                    B'_{t}(x) = \begin{cases}
                                            \eta^{t} \quad \text{if } x\in N(\bar{x}^{t}, k), \\
                                            0 \quad \text{otherwise}.
                                        \end{cases}
                    \label{eq:LB_Bounding_Function}
                \end{equation}
That is, $\eta^t$ is the optimal value of $\mathcal{S'}$ or, equivalently, the best upper bound accomplished by a solution in $N(\bar{x}^{t}, k)$. Although intuitively correct, the following must be proved. 
        \begin{lemma}
                    Bounding function (\ref{eq:LB_Bounding_Function}) satisfies properties $P.1$ and $P.2$.
                    \label{lem:1}
                \end{lemma}
        \proof{Proof}
        Let $N(\bar{x}, k) = \{x^{1}, x^{2},\ldots x^n\}$. Solving $\mathcal{S}'$ is equivalent to solving subproblem $\mathcal{S}$ for all $x^i, i=1,\ldots,n$. Equivalently, if $\eta^{t}$ is the optimal objective value of $\mathcal{S}'$, then $\eta^t = min_{i=1,\ldots,n}\{\zeta^i\geq f(x^i) + g(y))|$ $C_2(y), C_3(x^i, y), y\in D_y\}$. Hence, the objective value $\zeta^i$ is an upper bound of $\eta^{t}$ for $i=1,\ldots,n$ (i.e., for all $x \in N(\bar{x}^{t}, k)$).
        
        Consider a bounding function of the form (\ref{eq:Bounding_Function}) for each $x^i, i=1,\ldots,n$ and the associated upper bound $\zeta^i$. By definition, (\ref{eq:Bounding_Function}) adheres to properties $P.1$ and $P.2$, hence the cut (\ref{eq:benders_cut}) is valid for all solutions $x$ of $\mathcal{M}$. Notice that any solution $x^i$ for which $\zeta^i> \eta^{t}$ cannot be retrieved by the master problem $\mathcal{M}$ at any iteration following $t$, as some other solution has provided a smaller objective value $\eta^{t}$. 
        But then, constructing all $B_{t}(x^i), i=1,\ldots,n$ is equivalent in LBBD terms with constructing a single function (\ref{eq:LB_Bounding_Function}), hence the latter satisfying properties $P.1$ and $P.2$ implies that so does (\ref{eq:LB_Bounding_Function}).
        \endproof
        
Lemma \ref{lem:1} implies that solving subproblems $\mathcal{S}$ for \emph{all} $x\in N(\bar{x}^{t}, k)$ is necessary to sustain the compliance of (\ref{eq:LB_Bounding_Function}) to Properties $P.1$ and $P.2$ and could indeed be efficient if $N(\bar{x}^{t}, k)$ is small. As each solution of $\mathcal{S}$ implies the generation of the respective cut (\ref{eq:benders_cut}),  the exhaustive exploration of the neighbourhood triggers the construction of multiple cuts per integer solution (e.g., \cite{jeihoonian,rei09}). Although this accelerates convergence, meaning that a fewer number of solutions of $\mathcal{M}$ is required for optimality to be reached, cut generation should be monitored to avoid a severe slowdown of $\mathcal{M}$. Restricting the procedure so that only cuts which are expected to contribute on a significant improvement of the lower bound are generated, or solving $\mathcal{S}$ for a limited number of solutions of $N(\bar{\textsc{x}}^{t}, k)$ could be efficient alternatives but sacrifice optimality.
        
Our idea is that, as long as the best bound $\eta^t$ associated with the entire neighbourhood $N(\bar{x}, k)$ at solution $t$ is obtained, a single combinatorial cut capturing all solutions in the neighbourhood is sufficient, thus avoiding both the overload of adding one cut per solution and loss of optimality if choosing a fraction of them. Specifically, after solving $\mathcal{S}'$ for the neighbourhood $N(\bar{x}^{t}, k)$, thus obtaining the local bound $\eta^{t}$ (which will be the new incumbent one, if it dominates the previous one), then we add the local branching cut (\ref{eq:local_branch_cut}), so that any solution of the neighbourhood will be ignored during the subsequent exploration of the Branch-and-Cut tree.

%% file: scheduling.tex
\section{A decomposition scheme for machine scheduling} \label{section:problem}
    Let $\mathcal{P}$ be the problem where a set of unsplittable jobs $J$ are assigned and then sequenced to a set of parallel machines $M$. The machines are unrelated, hence there is a processing time $p_{jm}$ for job $j\in J$ on machine $m\in M$. The processing operation of a job should succeed a setup operation, hence the setup time of job $j$ is denoted by $s_{ijm}$ where $i\in J\setminus \{j\}$ indicates the preceding job on the same machine. Each machine handles at most one job at a time and no more than $R$ setup operations can be performed simultaneously.
    
    \subsection{Position-based MILP formulation} \label{sec:MILP}
        As each machine cannot execute more than one job at a time, $\mathcal{P}$ requires the so-called precedence constraints, so that the solution consists of valid sequences of jobs. Such constraints are typical `big-M' inequalities, which are known to imply weak lower bounds. This is why we opt for a \emph{position-based} MILP formulation, which allows us to obtain stronger lower bounds for the completion times of the assigned jobs.

        Specifically, we consider that each machine $m\in M$ has $|J|$ slots, where each slot can be occupied by at most one job. Binary variables $x_{ijm}$ are set to 1 if job $j\in J$ is assigned to slot $i$ of machine $m\in M$, 0 otherwise. The values of the binary variables determine the processing, setup and completion times of the respective slots, denoted by continuous variables. The annotations of the proposed MILP are explicitly presented in Table \ref{tab:milp_annotations}:

        \begin{table}[ht] \small
            \caption{Sets, parameters and variables}
            \centering
            \begin{tabular}[t]{ll}
                \hline
                \textbf{Sets} &\\
                \hline
                $J$         & Jobs \\
                $M$         & Machines \\
                $J^*$       & Jobs, except for the first one (used to indicate slots of machines) \\
                \hline
                \textbf{Parameters} &\\
                \hline
                $p_{jm}$    & Processing time of job $j\in J$ to machine $m\in M$ \\
                $s_{ijm}$   & Setup time of job $j\in J$ succeeding job $i\in J\setminus \{j\}$ to $m\in M$ \\ 
                $s^{-}_{jm}$ &A lower bound of $s_{ijm}$; $s^{-}_{jm} = min_{l\in J\setminus \{j\}}\{s_{jlm}\}$ \\
                $V$         & A large real number\\
                \hline
                \textbf{Variables} &\\
                \hline
                $x_{ijm}\in \{0,1\}$    & 1 if job $j\in J$ is assigned to slot $i\in J$ of machine $m\in M$, 0 otherwise \\
                $C_{im}\geq 0$          & Completion time of job or slot $i\in J$ to machine $m\in M$\\
                $P_{im}\geq 0$          & Processing time of slot $i\in J$ \\
                $S_{im}\geq 0$          & Setup time of slot $i\in J$ \\
                \hline
            \end{tabular}
            \label{tab:milp_annotations}
        \end{table}

For any regular objective $z$  the following MILP will be used as the master problem $\mathcal{M}$:
        
        \begin{flalign}
            \text{min } &z && \notag &&\\ 
            &\sum_{i\in J}\sum_{m\in M}x_{ijm} = 1 && \forall j\in J \label{eq:c1} &&\\
            &\sum_{j\in J}x_{ijm} \leq 1 && \forall i\in J, m\in M \label{eq:c2} &&\\
            &P_{im} = \sum_{j\in J}p_{jm}\cdot x_{ijm} && \forall i\in J, m\in M \label{eq:c3} &&\\
            &\sum_{k\in J\setminus \{j\}}s_{kjm}\cdot x_{i-1km} - S_{im}\leq V\cdot(1 - x_{ijm}) && \forall i\in J^{*}, m\in M, j\in J \label{eq:c4} &&\\
            &C_{im} = C_{i-1m} + P_{im} + S_{im} && \forall i\in J^{*}, m\in M \label{eq:c5} &&\\
            &C_{0m} = P_{0m} + S_{0m} && \forall m\in M \label{eq:c6} &&\\
            &\sum_{j\in J}x_{ijm} \geq \sum_{j\in J}x_{i-1jm} && \forall i\in J^{*}, m\in M \label{eq:c7} &&\\
            & && \notag &&\\
            &x_{ijm}\in \{0, 1\} && \forall i\in J, j\in J, m\in M \notag &&\\
            &P_{im}, S_{im}, C_{im}\geq 0 && \forall i\in J, m\in M \notag &&
        \end{flalign}

By (\ref{eq:c1}), each job is assigned to exactly one slot of one machine. Constraints (\ref{eq:c2}) ensure that each slot cannot be occupied by more than one job. Constraints (\ref{eq:c3}) determine the processing time of slot $i$ on machine $m$, as defined by the assigned job. As the setup times are sequence-dependent, Constraints (\ref{eq:c4}) ensure that, if job $j$ is assigned to a particular slot $i$ on machine $m$, then the setup time of the respective slot will be determined by the preceding job $k\neq j$. The completion time of each slot is obtained by the summation of the processing and setup times, added by the completion time of the preceding slot (Constraints (\ref{eq:c5}) and (\ref{eq:c6})). We note that subset $J^{*}$ denotes $J\setminus \{0\}$, and it is used to indicate the set of slots, except for the first one (otherwise, variables $x_{i-1km}$ or $C_{i-1m}$ could not be defined).
        
Constraints (\ref{eq:c7}) ensure the continuity of sequences, meaning that any slot cannot be occupied, as long as the succeeding one is empty. In other words, (\ref{eq:c7}) enforce that all empty slots appear first. Neglecting constraints (\ref{eq:c7}) could imply invalid values of completion times. E.g., if a slot $i$ remains vacant, while slot $i-1$ is occupied for any machine $m$, then the completion time $C_{im}$ would be equal to $C_{i-1m}$ by (\ref{eq:c5}), although its value should be normally set to 0. An indicative solution of $\mathcal{M}$ for a minimal instance of 5 jobs on 2 machines is shown in Figure \ref{fig:example}:

        \begin{figure}[H]
            \begin{center}\includegraphics{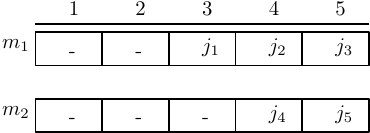}
                \caption{An indicative solution for an instance of $5$ jobs on $2$ machines}
                    \label{fig:example}
            \end{center}
        \end{figure}

\subsection{A family of valid inequalities}
As mentioned before, sequence-dependencies are formulated in Integer Programming as `big-M' constraints, which are known to provide weak dual bounds. As the setup times are sequence-dependent, such constraints are (\ref{eq:c4}). To obtain stronger lower bounds through our MILP formulation of $\mathcal{M}$, we consider a set of valid inequalities as follows.

Let $j\in J$ be a job assigned to slot $i-1\in J$ of machine $m\in M$ (i.e., $x_{i-1jm} = 1$). The assignment of $j$ to $i-1$ ensures that the candidate values of the setup time of the succeeding slot $i$ are from the set $\{s_{jkm}:k\in J\setminus\{j\}\}$. Given that, a straightforward lower bound of $S_{im}$ is the minimum value of this range: $S_{im} \geq min_{k\in J\setminus\{j\}}\{s_{jkm}\}$. Therefore, adding this bound as a linear inequality contributes on the enhancement of the lower bound of setup times (which are normally defined by the weak 'big-M' constraints), without loss of validity as no feasible solutions are eliminated. If $s^{-}_{jm} = min_{k\in J\setminus\{j\}}\{s_{jkm}\}$ $\forall j\in J, m\in M$, then:
        \begin{flalign}
            &S_{im}\geq \sum_{j\in J}s^{-}_{jm}\cdot x_{i-1jm} && \forall i\in J^{*}, m\in M \label{eq:c8} && &&
        \end{flalign}

Let us note that a second lower bound could be obtained, assuming that job $j\in J$ is assigned to slot $i$ on $m$, and the setup time of slot $i$ is greater than $min_{k\in J\setminus \{j\}}\{s_{kjm}\}$. As, however, the range of candidate values should also include the case of $j$ being the first job in the sequence, in which cases its setup time has been assumed to be $0$, this lower bound would have no effect on variables $S_{im}$. We still mention this as it could be useful if the first jobs had non-zero setup times.

 \subsection{Obtaining feasible schedules} \label{sec:subpro}
        The optimal solution $\bar{x}$ of $\mathcal{M}$ can be represented as a set of fixed sequences $\bar{M}(\bar{x})$, where each sequence $m\in \bar{M}$ is a tuple of job-to-slot assignments, where $m_i=j$ if and only if $\bar{x}_{ijm} = 1$. As mentioned before, the execution of setup tasks is subject to the availability of $R$ resources; setting up a job occupies one resource, which will be released after the completion of the setup task. After obtaining the sequences $\bar{M}$, as determined by a solution $\bar{x}$ of $\mathcal{M}$, we construct a Constraint Programming formulation to allocate resources in an optimal way.

        Specifically, we consider interval variables $\texttt{setup}_{m_{i}}$ and $\texttt{process}_{m_{i}}$, indicating the start and end times of the setup and processing of the $i^{th}$ job of machine $m$. For each machine, a sequence variable $\texttt{sequence}_{m}$ is constructed. To allocate resources, we consider the global constraint $\texttt{Cumulative}$. In details, the proposed subproblem $\mathcal{S}$ is formulated as it follows:

        \begin{flalign}
            \text{min } &\zeta && \notag &&\\ 
            &\texttt{noOverlap}(\texttt{sequence}_{m}) && \forall m\in \bar{M} \label{eq:s1} &&\\
            &\texttt{startAtEnd}(\texttt{process}_{m_{i}}, \texttt{setup}_{m_{i}}) && \forall m\in \bar{M}, i\in \{1, ..., |m|\} \label{eq:s2} &&\\
            &\texttt{previous}(\texttt{sequence}_{m}, \texttt{process}_{m_{i-1}}, \texttt{setup}_{m_{i}}) && \forall m\in \bar{M}, i\in \{2, ..., |m|\} \label{eq:s3} &&\\
            &\texttt{Cumulative}(\texttt{start}(\texttt{setup}_{m_{i}}), (\texttt{size}(\texttt{setup}_{m_{i}}), 1, R) && \label{eq:s4} &&\\
            & && \notag &&\\
            &\texttt{sequence}_{m}: \{\texttt{setup}_{m_{1}}, ..., \texttt{setup}_{m_{|m|}}, \texttt{process}_{m_{1}}, ..., \texttt{process}_{m_{|m|}}\} && \forall m\in \bar{M} \notag &&\\
            &\texttt{setup}_{m_{i}}, \texttt{process}_{m_{i}} : Interval && \forall m\in \bar{M}, i\in \{1, ..., |m|\} \notag &&
        \end{flalign}

The global constraint \texttt{noOverlap} (\ref{eq:s1}) ensures that setup and process tasks of the same sequence do not overlap. Constraints (\ref{eq:s2}) ensure that the processing of each job starts at the end of the respective setup, without delay, as the global constraint \texttt{startAtEnd} denotes. To the contrary, it is possible that the available resources do not allow a continuous succession between the completion of a job and the beginning of the setup of the following one; hence, we opt for the global constraint \texttt{previous} (\ref{eq:s3}). Finally, the allocation of resources is formulated through the global constraint \texttt{Cumulative} (\ref{eq:s4}), which forbids the simultaneous activation of more than $R$ setup tasks, each one occupying one unit of resources, which start at \texttt{start}(\texttt{setup}) and end after \texttt{size}(\texttt{setup}). We note that the size of interval variables are predefined by the respective processing and setup times, as the sequences of jobs on machines remain intact, as provided by the solution of $\mathcal{M}$.

Although the run-time of $\mathcal{S}$ depends on several parameters (e.g., the objective $\zeta$, the number of resources $R$, the size of setup times), it reaches optimality in no more than a few seconds in all experiments discussed in Section \ref{sec:exp}.

\subsection{The partitioning scheme}
Let $\mathcal{M}$ be the master problem of the decomposition scheme. We consider an integer solution, consisted of a vector of variables $\{\bar{x}^{t}_{ijm} = 1:i\in J, j\in J, m\in M\}$. After solving the respective subproblem $\mathcal{S}$, a new upper bound $\zeta^{t}$ is obtained. If $\zeta^{t}$ is better than the incumbent upper bound, then the vector of variables is stored as the incumbent solution, the value of $\zeta^{t}$ is set as the new incumbent bound, and the following cut is added to $\mathcal{M}$:
        \begin{flalign}
            \sum_{(i,j,m):\bar{x}^{t}_{ijm} = 1}(1 - x_{ijm}) \geq 1 \label{eq:f_cut} &&
        \end{flalign}
We observe that (\ref{eq:f_cut}) might be quite weak to imply a fast convergence, as only a single permutation of sequences of jobs is eliminated at each integer solution, which can easily be overrun by a slightly different one at an upcoming iteration, without a significant improvement of the lower bound. However, stronger optimality cuts, such as the ones shown by \cite{hooker07} or \cite{tran}, are no longer valid for $\mathcal{P}$. This is because machines are subject to the availability of the same resource, hence changing the sequence for one machine could affect the optimal schedule of others, even if their sequences remain intact. This is also observed by \cite{Pey20}, hence resource constraints are the ones making (\ref{eq:f_cut}) the sole candidate for a Benders cut.

On the positive side, Branch-and-Check remains valid, hence Algorithm 1, listed below, converges to optimality. Please note that each time $\mathcal{M}$ returns a new integer solution $\bar{x}^{t}$ at iteration $t$, it also returns the lower bound obtained through the LP-relaxations at all nodes of the Branch-and-Check tree enumerated up to that point.

        \begin{algorithm}[H]         
            \fontsize{9}{11}\selectfont
            \caption{A Branch-and-Check algorithm for $\mathcal{P}$}
            Let $LB = 0$ be the lower bound, $UB = \infty$ be the upper bound, $x'$ be the incumbent solution, and $t = 0$ be the number of integer solutions;\\
            While solving $\mathcal{M}$:
            \begin{algorithmic}[1]
                \WHILE {$LB < UB$}
                    \IF {a new integer solution $\bar{x}^{t}$ is obtained}
                        \STATE $LB \leftarrow$ lower bound of $\mathcal{M}$ once obtaining $\bar{x}^{t}$;
                        \STATE Let $\bar{z}^{t}$ be the objective value of $\mathcal{M}$ attained by $\bar{x}^{t}$;
                        \IF {$\bar{z}^{t} < UB$}
                            \STATE Solve $\mathcal{S}$ and let $\zeta^{t}$ be its objective value;
                            \IF {$\zeta^{t}<UB$}
                                \STATE $UB \leftarrow \zeta^{t}$,  $x'\leftarrow \bar{x}^{t}$;
                            \ENDIF
                            \STATE Add (\ref{eq:f_cut}) to $\mathcal{M}$;
                            \STATE Set $t = t + 1$;
                        \ENDIF
                    \ENDIF
                \ENDWHILE
                \STATE Return $LB$ and $x'$;
            \end{algorithmic}
            \label{alg:regular_BnC}
        \end{algorithm}

%% file: local_branching.tex
\section{Strengthening the decomposition} \label{sec:lb}
    \subsection{Locally Branch and Check}
        The \textit{swap} operator is quite standard in local search algorithms hence also called a swap \textit{movement}. Given a schedule, a swap exchanges only two jobs assigned to the same machine, leaving all other jobs intact. A sequence of $n$ jobs in a single machine allows for $\frac{n(n-1)}{2}$ swaps, i.e., the number of combinations of $2$ out of $n$ objects. Although assigning the jobs to multiple machines blurs that calculation, one can easily prove that the number of swaps is still no larger than $\frac{n(n-1)}{2}$ hence the swap neighbourhood remains $O(n^2)$. The literature calls this movement an \textit{internal swap} to emphasise that swaps are allowed only within a machine, i.e., an \textit{external swap} would switch two jobs assigned to different machines.

        Let $\bar{x}^t$ be the optimal solution of $\mathcal{M}$ at iteration $t$, called also the \emph{reference solution} of $t$.  The \emph{local branching constraint} for the $k$-OPT neighbourhood of that solution, i.e., constraint (\ref{eq:local_branch_cut}) has the form:
        \begin{equation}
            \sum_{m\in M}\sum_{(i,j):\bar{x}^{t}_{ijm} = 1}(1 - x_{ijm}) + \sum_{m\in M}\sum_{(i,j):\bar{x}^{t}_{ijm} = 0}x_{ijm} \leq k \label{eq:k_constraint}
        \end{equation}

    The structure of $\mathcal{P}$ imposes a constant sum of all binary variables, since all $|J|$ jobs are assigned to exactly one slot of one machine. In more detail, adding up constraints (\ref{eq:c1}) for all $j \in J$ yields:
         \begin{equation}
            \sum_{m\in M}\sum_{(i,j):\bar{x}^{t}_{ijm} = 1}x_{ijm} + \sum_{m\in M}\sum_{(i,j):\bar{x}^{t}_{ijm} = 0}x_{ijm} = |J| \label{eq:constant}
         \end{equation}
Therefore, (\ref{eq:k_constraint}) can be simplified (as in \cite{fischetti03}):
         \begin{equation}
            \sum_{m\in M}\sum_{(i,j):\bar{x}^{t}_{ijm} = 1}(1 - x_{ijm}) \leq \frac{k}{2} \label{eq:k'_constraint}
         \end{equation}

To employ (\ref{eq:k'_constraint}) within Branch-and-Check requires identifying a neighbourhood that is $k$-OPT for some value of $k$. Notice that such a neighbourhood is not only problem-specific but also formulation-specific. We need a small change in our formulation in order to define a $4$-OPT neighbourhood that includes internal swaps plus a few other movements.


Let $j$ and $j'$ be two jobs assigned to the same machine $m$ at any integer solution $\bar{x}^{t}$ (i.e., $\bar{x}^{t}_{ijm} = 1$, $\bar{x}^{t}_{i'j'm} = 1$). After implementing an internal swap operator, we notice that variables $x_{ijm}$ and $x_{i'j'm}$ are switched from 1 to 0, and variables $x_{i'jm}$ and $x_{ij'm}$ are switched from 0 to 1, meaning that the internal swap movements could be a part of a $4$-OPT neighbourhood for a subset of variables $x$. Notice however that the same holds also for external swaps, thus making the $4$-OPT neighbourhood rather large for our purposes.


To exclude external swaps, let us append to $\mathcal{M}$ the vector of auxiliary variables $y$ as follows: 
\begin{equation} \label{eq:yvars}
    y_{jm} = \sum_{i\in J}x_{ijm}\text{,  } \forall j \in J, m \in M,
\end{equation}
\begin{equation*}
    y_{jm} \in \{0,1\} \text{,  } \forall j \in J, m \in M.
\end{equation*} 
\noindent That is, variable $y_{jm}$ indicate whether job $j$ is assigned to machine $m$. 

Let $\mathcal{M}'$ be the formulation of $\mathcal{M}$, as listed in Section \ref{sec:MILP}, extended by the vector of variables $y$ and equalities (\ref{eq:yvars}). The local branching constraint (\ref{eq:k'_constraint}) for formulation $\mathcal{M}'$ becomes

\begin{equation}
            \sum_{m\in M}\sum_{j:\bar{y}^{t}_{jm}=1}(1 - y_{jm}) + \sum_{m\in 
 M}\sum_{(i,j):\bar{x}^{t}_{ijm} = 1}(1 - x_{ijm}) \leq \frac{k}{2}. \label{eq:local_branching_constraint}
\end{equation}

Notice that the value of any variable $y_{jm}$ is invariant to internal swaps. Therefore, 
an internal swap remains inside the $4$-OPT neighbourhood for a given solution $(\bar{x}^{t}, \bar{y}^{t})$. As mentioned above, any internal swap will imply switching of two variables from $0$ to $1$ and two different variables from $1$ to $0$, all four of them belonging in $\bar{x}^{t}$. At the same time, the values of variables $y$ will remain intact, as all jobs are assigned to the same machines. To the contrary, an external swap would still imply four switches of variables in $\bar{x}^{t}$, but it would also switch values of variables in $\bar{y}^{t}$, as different assignments of jobs to machines would occur, meaning that the defined $4$-OPT neighbourhood excludes them. Let us just emphasize that, although the `trick' of introducing vector $y$ makes $4$-OPT smaller, one can easily revert to $\mathcal{M}$ and use the local branching idea within Branch-and-Check but for a much larger $4$-OPT.  

\begin{figure}[H]
            \begin{center}\includegraphics{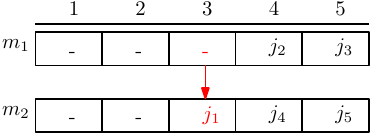}
                \caption{An indicative example of moving the starting job $j_{1}$ of machine $m_{1}$ to the starting slot of machine $m_{2}$}
                    \label{fig:single}
            \end{center}
        \end{figure}

Still, to ensure the validity of the $4$-OPT neighbourhood for formulation $\mathcal{M}'$, we should consider all candidate solutions which have a symmetric difference of at most $4$ variables in $(\bar{x}^{t}, \bar{y}^{t})$. For that purpose, we define another movement called \emph{starting-job shift} (in lack of a literature term, to our knowledge at least). This movement means that the starting job of a machine $m$ is shifted to become the starting job of a machine $m'\neq m$, as shown indicatively in Figure \ref{fig:single}.
        
\begin{definition} \label{def:neighbour}
            The neighbourhood $N(\bar{x}^{t}, \bar{y}^{t})$ includes \textbf{(i)} the $t^{th}$ solution $(\bar{x}^{t}, \bar{y}^{t})$ of $\mathcal{M'}$, \textbf{(ii)} solutions obtained by internal swaps for all pairs of jobs $j$ and $j'$ assigned to the same machine $m$ in $(\bar{x}^{t}, \bar{y}^{t})$, and \textbf{(iii)} solutions obtained from $(\bar{x}^{t}, \bar{y}^{t})$ by starting job-shifts for all pairs of machines $m$ and $m'$.
    \end{definition}

\begin{proposition}
            $N(\bar{x}^{t}, \bar{y}^{t})$ coincides with the $4$-OPT neighbourhood for $\mathcal{P}$. \label{prop:1}
\end{proposition}  
 \proof
We prove that a solution in $N(\bar{x}^t, \bar{y}^t)$ satisfies (\ref{eq:local_branching_constraint}) and any other solution violates it. 

Looking at Definition \ref{def:neighbour}(i), the reference solution $(\bar{x}^t, \bar{y}^t)$ makes the left-hand side of (\ref{eq:local_branching_constraint}) $0$. To see that, notice that $\sum_{m\in M}\sum_{j:\bar{y}^{t}_{jm} = 1} \bar{y}_{jm} = |J|$ as all jobs are assigned to machines, hence $\sum_{m\in M}\sum_{j:\bar{y}^{t}_{jm} = 1}(1 - \bar{y}_{jm}) = 0$ and, similarly, $\sum_{m\in M}\sum_{(i,j):\bar{x}^{t}_{ijm} = 1}\bar{x}_{ijm} = |J|$ thus the second term in (\ref{eq:local_branching_constraint}) becomes also $0$.

Having discussed above that internal swaps switch exactly two variables from $0$ to $1$ shows that all solutions complying to Definition \ref{def:neighbour}(ii) satisfy (\ref{eq:local_branching_constraint}) at equality. It remains to check starting-job shifts. Such a movement shifts the job $j$ occupying the first slot of machine $m$ to the first slot of machine $m'$ hence changing variable $y_{jm}$ from $1$ to $0$ and variable $y_{jm'}$ from $0$ to $1,$ the same occurring for variables $x_{1jm}$ and $x_{1jm'}.$ No other variable changes its value, as  (\ref{eq:c7}) enforces that all empty slots  appear first. Therefore, for any solution complying with Definition \ref{def:neighbour}(iii) it holds that $\sum_{m\in M}\sum_{j:\bar{y}^{t}_{jm}=1}y_{jm} = \sum_{m\in M}\sum_{(i,j):\bar{x}^{t}_{ijm} = 1}x_{ijm} = |J|-1$; hence both terms in left-hand side of (\ref{eq:local_branching_constraint}) receive a value $1$ and (\ref{eq:local_branching_constraint}) is satisfied at equality.

Let us now show that any solution not in $N(\bar{x}^t, \bar{y}^t)$ violates (\ref{eq:local_branching_constraint}). Let $(x', y')\notin N(\bar{x}^t, \bar{y}^t)$ be a feasible solution of $\mathcal{M}$. As $x'\neq \bar{x}^t$, there is some triple $(i,j,m)$ such that $x'_{ijm} = 0$ whereas $\bar{x}^{t}_{ijm} = 1$. In other words, in solution $(x', y')$, some job $j$ is definitely transferred to another slot $i'$ of a (possibly different from $m$) machine $m'$. The implication up to here is that $$\sum_{m\in M}\sum_{(i,j):\bar{x}^{t}_{ijm} = 1}(1 - x'_{ijm}) \geq 1 \text{ for all } (x', y')\neq (\bar{x}^t, \bar{y}^t),$$ 
\noindent therefore, for (\ref{eq:local_branching_constraint}) not to be violated, we should have 
$$\sum_{m\in M}\sum_{j:\bar{y}^{t}_{jm}=1}(1 - y'_{jm}) \leq 1.$$
\noindent We examine two cases, distinguished by whether some job is transferred to a different machine.

If at least one job $j$ is transferred from $m$ to a machine $m'\neq m$, then $y'_{jm}=0$ whereas $\bar{y}^{t}_{jm} = 1$, yielding that $\sum_{m\in M}\sum_{j:\bar{y}^{t}_{jm}=1}(1 - y'_{jm}) \geq 1$. This, in combination with the previous inequality, yields $\sum_{m\in M}\sum_{j:\bar{y}^{t}_{jm}=1}(1 - y'_{jm}) = 1,$ which in turn implies that all other jobs $j'\neq j$ remain in the same machine in both solutions. All these jobs (i.e., all but $j$) must also remain in the same slots (i.e., $x'$ coincides with $\bar{x}$ except for entries $(i,j,m)$ and $(i',j,m')$), as otherwise $\sum_{m\in M}\sum_{(i,j):\bar{x}^{t}_{ijm} = 1}(1 - x'_{ijm})$ would be at least $2,$ thus violating (\ref{eq:local_branching_constraint}). But then,  Constraints (\ref{eq:c7}) imply that $j$ occupies the first slot in machine $m'$ (as empty slots appear first and not between occupied ones), thus $(x',y')$ is obtainable from $(\bar{x}^t, \bar{y}^t)$ through a starting-job shift, contradicting that $(x',y')\notin N(\bar{x}^t, \bar{y}^t)$.

\begin{figure}[H]
            \begin{center}\includegraphics{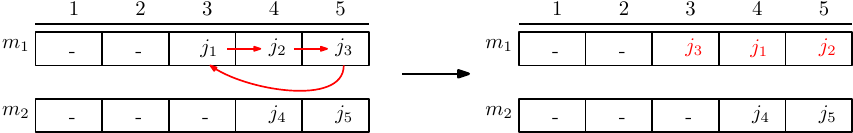}
                \caption{A slot change by three jobs}
                    \label{fig:non_swap}
            \end{center}
        \end{figure}

Otherwise, all jobs remain in the same machine in both solutions, i.e.,  $\bar{y}= y'$, hence $\sum_{m\in M}\sum_{j:\bar{y}^{t}_{jm}=1}(1 - y_{jm}) =0.$ This means that job $j$ is transferred to a different slot $i'$ of the same machine $m.$ Then, some other jobs in machine $m$ must also change their slots, as no intermediate vacancies can appear by Constraints (\ref{eq:c7}). Should two or more jobs other than $j$ change their slots in machine $m$, it would yield $\sum_{m\in M}\sum_{j:\bar{y}^{t}_{jm}=1}(1 - y_{jm}) +  \sum_{m\in M}\sum_{(i,j):\bar{x}^{t}_{ijm} = 1}(1 - x_{ijm})\geq 3$, in violation of (\ref{eq:local_branching_constraint}) (see a minimal example in Figure \ref{fig:non_swap}).
\noindent For only one job other than $j$ to change its slot without leaving an empty slot, $j$ should move into a previously occupied slot $i'$, forcing job $j'$ for which $\bar{x}^{t}_{i'j'm} = 1$ to move into slot $i$. But then, $(x',y')$ is obtainable from $(\bar{x},\bar{y})$ by an internal swap, contradicting $(x', y')\notin N(\bar{x}^t, \bar{y}^t)$. 
Overall, any solution $(x', y')\notin N(\bar{x}^t, \bar{y}^t)$ violates (\ref{eq:local_branching_constraint}) and the proof is complete.
Let us now show that any solution not in $N(\bar{x}^t, \bar{y}^t)$ violates (\ref{eq:local_branching_constraint}). Both terms 
$\sum_{m\in M}\sum_{j:\bar{y}^{t}_{jm}=1}y_{jm}$ and $\sum_{m\in M}\sum_{(i,j):\bar{x}^{t}_{ijm} = 1}x_{ijm}$ receive value $|J|$, $|J|-1$ and $|J|-2$ for a solution $(x,y)$ complying with Definition \ref{def:neighbour}(i), (ii) and (iii) respectively. Also, any solution having value at most $|J|-3$ for either term yields one of the terms in the left-hand side of (\ref{eq:local_branching_constraint}) at value $3$ hence violating that inequality. As the term $\sum_{m\in M}\sum_{(i,j):\bar{x}^{t}_{ijm} = 1}x_{ijm} = |J|$ is valid only for the reference solution, let us examine only solutions for which 
the term $\sum_{m\in M}\sum_{(i,j):\bar{x}^{t}_{ijm} = 1}x_{ijm}$ becomes $|J|-1$  or $|J|-2.$
\endproof
Note that replacing variables $y_{jm}$ in (\ref{eq:local_branching_constraint}) through (\ref{eq:yvars}), projecting out vector $y$ in polyhedral terms, yields an equivalent local branching constraint for $\mathcal{M}$ using only variables $x$.

\subsection{Exploring faster through domination} \label{section:dom}
Let us now describe how we integrate the local branching framework into the Branch-and-Check Algorithm \ref{alg:regular_BnC}. Each time a new integer solution of $\mathcal{M}$ is found, this solution is used as a \emph{reference solution} to construct the $4$-OPT neighbourhood of Definition \ref{def:neighbour}. To retrieve the best upper bound $\eta$, all solutions of the neighbourhood must be explored, meaning that the subproblem $\mathcal{S}$ should be solved for all of them,. Although solving $\mathcal{S}$, as listed in of Section \ref{sec:subpro}, requires no more than a few seconds, even the number of internal swaps can easily become more than one thousand. 

To speed up this exploration, given a solution $(x,y)\in N(\bar{x}^t, \bar{y}^t)$, we could drop resource constraints and obtain the `resource-free' relaxation of $\mathcal{S}$. This would compute a lower (i.e., an optimistic) bound $\zeta'$ on the value $\zeta$ of the subproblem for that particular solution. If that is larger than the best upper bound (UB) found among the subproblems solved so far, there is no reason to solve $\mathcal{S}$ for that solution $(x,y).$

We also opt for a CP formulation that captures all candidate solutions of the neighbourhood, instead of solving the subproblem $\mathcal{S}$ of Section \ref{sec:subpro} for each one of them. For each $s\in N(x^{t}, y^{t})$, let $\bar{M}_{s}$ be the set of the respective sequences of jobs. Consider the following formulation, referred to hereafter as $\mathcal{S}'$.

    \begingroup
        \footnotesize \selectfont
        \begin{flalign}
            \text{min } &\sum_{s\in N(x^{t}, y^{t})}\zeta_{s}\cdot \lambda_{s} && \notag &&\\
            &\sum_{s\in N(\bar{x}^t, \bar{y}^t)}\lambda_{s} = 1 && \label{eq:t1} &&\\
            &\texttt{noOverlap}(\texttt{sequence}_{sm}) && \forall s\in N(\bar{x}^t, \bar{y}^t), m\in \bar{M}_{s} \label{eq:t2} &&\\
            &\texttt{startAtEnd}(\texttt{process}_{sm_{i}}, \texttt{setup}_{sm_{i}}) && \forall s\in N(\bar{x}^t, \bar{y}^t), m\in \bar{M}_{s}, i\in \{1, ..., |m|\} \label{eq:t3} &&\\
            &\texttt{previous}(\texttt{sequence}_{sm}, \texttt{process}_{sm_{i-1}}, \texttt{setup}_{sm_{i}}) && \forall s\in N(\bar{x}^t, \bar{y}^t), m\in \bar{M}_{s}, i\in \{2, ..., |m|\} \label{eq:t4} &&\\
            &\texttt{Cumulative}(\texttt{start}(\texttt{setup}_{sm_{i}}), (\texttt{size}(\texttt{setup}_{sm_{i}}), 1, R) && \forall s\in N(\bar{x}^t, \bar{y}^t) \label{eq:t5} &&\\
            &\texttt{presenceOf}(\texttt{setup}_{sm_{i}}) = \lambda_{s} && \forall s\in N(\bar{x}^t, \bar{y}^t), m\in \bar{M}_{s}, i = \{1, ..., |M|\} \label{eq:t6} &&\\
            &\texttt{presenceOf}(\texttt{process}_{sm_{i}}) = \lambda_{s} && \forall s\in N(\bar{x}^t, \bar{y}^t), m\in \bar{M}_{s}, i = \{1, ..., |M|\} \label{eq:t7} &&\\
            &\zeta_{s} \geq UB \rightarrow \lambda_{s} = 0 && \forall s\in N(\bar{x}^t, \bar{y}^t) \label{eq:t8} &&\\
            & && \notag &&\\
            &\texttt{sequence}_{sm}: \{\texttt{setup}_{sm_{1}}, ..., \texttt{setup}_{sm_{|m|}}, \texttt{process}_{sm_{1}}, ..., \texttt{process}_{sm_{|m|}}\} && \forall s\in N(\bar{x}^t, \bar{y}^t), m\in \bar{M}_{s} \notag &&\\
            &\texttt{setup}_{sm_{i}}, \texttt{process}_{sm_{i}} : Interval && \forall s\in N(\bar{x}^t, \bar{y}^t), m\in \bar{M}_{s}, i\in \{1, ..., |m|\} \notag &&
        \end{flalign}
        \endgroup

The above formulation $\mathcal{S}'$ considers a set of binary variables $\lambda_{s}$, indicating whether the solution $s\in N(\bar{x}^t, \bar{y}^t)$ is the locally optimal one. All interval variables are set as \emph{optional}, meaning that they will receive values of start and end times only if the respective binary variable is set to $1$ by (\ref{eq:t6})-(\ref{eq:t7}). We have also added a set of cut-off constraints (\ref{eq:t8}), ensuring that no solutions of higher upper bounds than the incumbent one (i.e., $UB$) will be selected.

Although this modified formulation is significantly larger than the one of Section \ref{sec:subpro}, solving it once is sufficient to obtain the best upper bound in the neighbourhood. Hence, it allows us to avoid an exhaustive exploration of the neighbourhood, which is more time-consuming for instances of practical scale. We also note that $\mathcal{S}'$ could be infeasible if no solution provides an improved upper bound, as (\ref{eq:t8}) would then contradict (\ref{eq:t1}).

\subsection{Improving Branch-and-Check with Local Branching}
Given a solution $\bar{x}^{t}$, local Branching creates two branches, one featured with (\ref{eq:k'_constraint}) and the second with $$\sum_{m\in M}\sum_{(i,j):\bar{x}^{t}_{ijm} = 1}(1 - x_{ijm}) \geq \frac{k}{2} + 1.$$ Our approach adds the latter inequality as a Benders cut, thus making $\mathcal{M}$ explore the second branch through Branch-and-Check. The first branch is explicitly enumerated through local search and then by solving $\mathcal{S}'$. 
That is, after obtaining the local optimum of the branch by solving $\mathcal{S}'$, we can safely prune it, so that the exploration of the master problem will search for solutions of a symmetric distance of at least $k+1$ from  $\bar{x}^{t}$. Specifically, after the implicit enumeration of the $4$-OPT neighbourhood of Definition \ref{def:neighbour}, we add the  cut:
        \begin{equation}
            \sum_{m\in M}\sum_{j:\bar{y}^{t}_{jm}=1}(1 - y_{jm}) + \sum_{m\in 
            M}\sum_{(i,j):\bar{x}^{t}_{ijm} = 1}(1 - x_{ijm}) \geq 3, \label{eq:lb_cut}
        \end{equation}
which again can be written in terms of only $x$ by replacing $y$ through (\ref{eq:yvars}). Although the approach discussed in Section \ref{section:dom} induces a faster exploration of $4$-OPT neighbourhoods, eventually the repetitive construction of such neighbourhoods for all integer solutions of $\mathcal{M}$ might become more burdensome than the Branch-and-Check scheme of Algorithm 1. Therefore, we opt for a more conservative implementation of  local branching only whenever this is more likely to accelerate LBBD convergence, i.e., only for integer solutions that fail to improve the incumbent upper bound. This enhanced Branch-and-Check scheme is shown below as Algorithm \ref{alg:lb_BnC}:

        \begin{algorithm}[H]         
            \fontsize{9}{11}\selectfont
            \caption{A Branch-and-Check algorithm with Local Branching for $\mathcal{P}$}
            Let $LB = 0$ be the lower bound, $UB = \infty$ be the upper bound, $x'$ be the incumbent solution, and $t = 0$ be the number of integer solutions;\\
            While solving $\mathcal{M}$:
            \begin{algorithmic}[1]
                \WHILE {$LB < UB$}
                    \IF {a new integer solution $(\bar{x}^t, \bar{y}^t)$ is obtained}
                        \STATE $LB \leftarrow$ lower bound of $\mathcal{M}$ once obtaining $\bar{x}^{t}$;
                        \STATE Let $\bar{z}^{t}$ be the objective value for solution $(\bar{x}^t, \bar{y}^t)$
                        \IF {$\bar{z}^{t} < UB$}
                            \STATE Solve $\mathcal{S}$ for solution $(\bar{x}^t, \bar{y}^t)$ and let $\zeta^{t}$ be its objective value;
                            \IF {$\zeta^{t}<UB$}
                                \STATE $UB \leftarrow \zeta^{t}$,  $x'\leftarrow (x^{t}, y^{t})$;
                                \STATE Add (\ref{eq:f_cut}) to $\mathcal{M}$;
                            \ELSE
                                \STATE Construct the neighbourhood of Definition \ref{def:neighbour};
                                \STATE Solve $\mathcal{S}'$ and let $\zeta^{t}$ be its objective value;
                                \IF {$\zeta^{t}(N) < UB$}
                                    \STATE $UB \leftarrow \zeta^{t}$,  $x'\leftarrow s: \lambda_{s} = 1$ from $\mathcal{S}'$;
                                \ENDIF
                                \STATE Add (\ref{eq:lb_cut}) to $\mathcal{M}$;
                            \ENDIF
                            \STATE Set $t = t + 1$;
                        \ENDIF
                    \ENDIF
                \ENDWHILE
                \STATE Return $LB$ and $x'$;
            \end{algorithmic}
            \label{alg:lb_BnC}
        \end{algorithm}

\noindent Observe that, as long as the integer solutions of $\mathcal{M}$ imply an improved upper bound, we construct regular cuts (\ref{eq:f_cut}). To the contrary, if a newly obtained solution does not improve the upper bound, we implement the local search operators of Definition \ref{def:neighbour} to search for an improvement in a neighbouring solution. After the exploration of the neighbourhood, we proceed to the generation of a local branching cut (\ref{eq:lb_cut}).

 \subsection{Warm-start heuristics} \label{section:heuristics}
We use three heuristics to provide a starting solution to the master problem $\mathcal{M}$, hence neglecting the resource constraints. The first heuristic, called GH-EDD, assigns jobs on machines using the MIP formulation (called \textit{MIP 3}) from \cite{ijpr2022}, and then sequences the assigned jobs on each machine based on the Earliest Due Date (EDD) rule. 
        
The second heuristic, named GH-slack, uses the same MIP as GH-EDD for the job assignment, but a different sequencing rule. Specifically, sequencing occurs from repeatedly calculating the formula $I_{jm} =  d_j - (load_m + p_{jm} + s_{ijm})$, where $d_j$ is the due date of job j, $load_m$ is the current load of machine $m$, $p_{jm}$ is the processing time of job $j$ on machine $m$ and $s_{ijm}$ is the setup time of $j$ on machine $m$, job $i$ being the one just sequenced on $m$. At each iteration, GH-slack  chooses the job $j$ with the minimum value $I_{jm}$ (i.e., the minimum `distance' from its deue date) and assigns it as the next job to $m$. 
        
The third heuristic imitates the well-known process described by \cite{Chen09}, which in turn extends the Apparent Tardiness Cost with Setups (ATCS) algorithm \cite{Pinedo97} for the case of unrelated parallel machines. The solution obtaining the lowest objective value (either the sum of completion times or the sum of tardiness) is fetched to the master problem.

%% file: experiments.tex
\section{Numerical Evidence} \label{sec:exp}
       \subsection{Two objective functions}
We consider two well-known objectives, which are the sum of completion times and total tardiness. For both objectives, we are interested in the efficiency of the proposed LBBD and Branch-and-Check (Algorithm \ref{alg:regular_BnC}) and the achieved by local branching enhancements listed in Section \ref{sec:lb} (Algorithm \ref{alg:lb_BnC}).

\noindent \textbf{Sum of completion times. }For the sum of completion times, the formulation of $\mathcal{M}$ includes an adequate number of variables and constraints. The auxiliary variable $z$, which denotes the objective of the master problem, is set to $$z = \sum_{m\in M}\sum_{i\in J}C_{im}.$$ For the subproblem $\mathcal{S}$ of Section \ref{sec:subpro}, variable $\zeta$ is set to $$\zeta = \sum_{m\in \bar{M}}\sum_{i\in \{1, ..., |m|\}}\texttt{endOf}(\texttt{process}_{m_{i}}).$$ For the subproblem $\mathcal{S}'$ of Section \ref{section:dom}, each variable $\zeta_{s}$ is set to $$\zeta_{s} =  \sum_{m\in \bar{M}_{s}}\sum_{i\in \{1, ..., |m|\}}\texttt{endOf}(\texttt{process}_{sm_{i}}) \text{,  } \forall s\in N(\bar{x}^{t}, \bar{y}^{t}).$$

\noindent \textbf{Sum of tardiness. }Several adjustments are required here as well. As the value of completion times of jobs are assigned to the variables $C_{im}$ of the respective slots $i\in J$, the values of due-date and tardiness should also be indicated by similar variables. We introduce the continuous variables $D_{im}\geq 0$ that indicate the due-date of slot $i\in J$ on machine $m\in M$. Their values are determined by the following constraints:
        \begin{equation}
            D_{im} = \sum_{j\in J}d_{j}\cdot x_{ijm} \qquad \forall i\in J, m\in M, \label{eq:tard1}
        \end{equation}
        where $d_{j}$ denotes the due-date of job $j\in J$. To define the value of tardiness of each slot, let $T_{im}\geq 0$ be a set of non-negative, continuous variables for all slots $i\in J$ on machines $m\in M$:
        \begin{equation}
            T_{im} \geq C_{im} - D_{im} \qquad \forall i\in J, m\in M \label{eq:tard2}
        \end{equation}
        The objective is the sum of variables $T_{im}$: $z = \sum_{m\in M}\sum_{i\in J}T_{im}$. For the subproblem $\mathcal{S}$, let $T^{*}_{j}$ be integer variables which indicate the value of tardiness of job $j\in J$. We also add the following constraints:
        \begin{equation}
            T^{*}_{m_{i}} = max\{0, \texttt{endOf}(\texttt{process}_{m_{i}})\} \qquad \forall m\in \bar{M}, i\in \{1, ..., |m|\} \label{eq:tard3}
        \end{equation}
        The objective function is set to $\zeta = \sum_{j\in J}T^{*}_{j}$.

\subsection{Generation of instances}
We have thoroughly considered which among an enormous number of possible instances should serve both as a representative set guaranteeing the validity of our conclusions and as a varying set showing the versatility of our approach. Under this dual aim, we have generated random instances with the following combinations of number of jobs $|J|$ and machines $|M|$: $|J| = 50 \text{ on } |M| = 5, |J| = 50 \text{ on } |M| = 10$, and $|J| = 100 \text{ on } |M| = 10$. For the number of resources $R$, we assume the values $2/5$ and $3/5$ of the number of machines, i.e., an instance with $|M| = 5$ is solved for $R \in \{2,3\}$ and an instance with $|M| = 10$ is solved for $R \in \{4,6\}$.
    
Regarding the processing times $p_{jm}$ of job $j\in J$ on machine $m\in M$, we consider the formula $b_j \cdot a_{jm} + \epsilon$, where $b_j$ and $a_{jm}$ are selected uniformly at random (u.a.r.) from $[1, 10]$, as in \cite{Fotakis16}, then increased by a `noise' value $\epsilon$, selected u.a.r. from $[0, 10]$. 
We consider three generators of setup times, denoted by a value $\alpha \in\{ 0, 1, 2\}$. For $\alpha = \{0, 1\}$, inspired by a similar generator \cite{Lee21}, $s_{ijm}$ are determined as a fraction of the respective value of $p_{jm}$: $s_{ijm} = \beta_{ijm}\cdot p_{jm}$ where $\beta_{ijm}$ is u.a.r. in [0.1, 0.5] and in [0.5, 1.0] for $\alpha = 0$ and $1,$ respectively. Regarding $\alpha = 2$, the values of setup times are \emph{independent} of the respective values of processing times, hence $s_{ijm}$ is u.a.r. in $[5,25]$ as in \cite{Bek19}.
                
To determine the values of due-date $d_{j}$ of jobs $j\in J$, we consider the $\tau$ and $\rho$ values of \cite{Dog79}, which are regularly used for generating due-date related instances in the literature. In particular, all instances assume $\rho = 0.2$ and $\tau \in \{0.5, 0.8\},$ the former value denoting more loose due-dates and the latter stricter ones. For each job $j\in J$, the value of due-dates is set u.a.r. from $\big( C_{\max}\cdot (1 - \tau - \frac{\rho}{2})$, $C_{\max}\cdot (1 - \tau + \frac{\rho}{2})\big)$, in which $C_{\max} = \frac{\sum_{j}({P^1_{j}})+ S^1}{|M|}$. $P^{1}_{j}$ is equal with $\min_{m\in M} \{p_{jm}\}$ and $S^{1}$ is set to $\sum_{j\in J^{\prime}}\min_{i\in J\setminus \{j\}, m\in M}\{s_{ijm}\}$, $J^{\prime}$ being the subset of the first $|J| - |M|$ jobs in ascending order. This is due to the assumption that the starting jobs require no setup, therefore the last $|M|$ jobs of $J$ in ascending order of setup times are not involved in the computation of an estimation $S^{1}$ for the average estimated makespan $C_{\max}$.

Table \ref{tab:generator} summarizes the values of all parameters for the randomly generated instances. Instances are denoted as $J\_M\_\tau \_\alpha$, indicating the number of jobs, the number of machines, the due-date tightness value $\tau$ and the setup parameter $\alpha$. For example, the instance $J50\_M5\_\tau0.5\_\alpha0$ consists of 50 jobs and 5 machines, has a due-date of tightness $\tau = 0.5$, and belongs to the group of $\alpha = 0$ setup times.
        
        \begin{table}[ht]
            \caption{Parameters of the Generator of Instances}
            \centering
            \begin{tabular}[t]{ll}
                \hline
                \textbf{Annotation} &\textbf{Uniform}\\
                \hline
                $|J|$               & $\{50, 100\}$ \\
                $|M|$               & $\{5, 10\}$ \\
                $R$                 & $\{\frac{2}{5}\cdot |M|, \frac{3}{5}\cdot |M|\}$ \\
                $a_{jm}$            & $U(0, 10)$ \\
                $b_{j}$             & $U(0, 10)$ \\
                $\epsilon$          & $U(0, 10)$ \\
                $p_{jm}$            & $a_{jm}\cdot b_{jm} + \epsilon$ \\
                $\alpha$            & $\{0, 1, 2\}$ \\
                $\beta_{ijm}$           & $U(0.1, 0.5)\cdot p_{jm}$ if $\alpha = 0$; $U(0.5, 1.0)\cdot p_{jm}$ if $\alpha = 1$\\
                $s_{ijm}$           & $\beta_{ijm}\cdot p_{jm}$ if $\alpha \in \{0, 1\}$; $U(5, 25)$ if $\alpha = 2$\\
                $\tau$              & $\{0.5, 0.8\}$\\
                $\rho$              & 0.2\\
                $P_{j}^{1}$         & $\text{min}_{j\in J} \{p_{jm}\}$\\
                $S^{1}$             & $\sum_{j\in J^{|J|-|M|}}\min_{i\in J\setminus \{j\}}\{s_{ijm}\}$ \\
                $C_{\max}$           & $\frac{\sum_{j}({P^1_{j}})+ S^1}{|M|}$. $P^{1}_{j}$ \\
                $d_{j}$             & $U(C_{\max}\cdot (1 - \tau - \frac{\rho}{2})$, $C_{\max}\cdot (1 - \tau + \frac{\rho}{2}))$\\
                \hline
            \end{tabular}
            \label{tab:generator}
        \end{table}

    \subsection{Results}

All experiments have been performed on a server with 4 Intel(R) Xeon(R) E-2126G @ 3.30GHz processors and 11 GB RAM, running CentOS/Linux 7.0. To solve the MILP-formulated master problem, we use the \emph{Gurobi 10.0.0} optimizer, summoned by the open-source optimisation library \emph{Pyomo 5.7.3} on \emph{Python 3.8.5}. The CP-formulated subproblem is solved by the \emph{CP Optimizer} of the \emph{DOCplex} module. A time limit of 3600 seconds is imposed on both Algorithms \ref{alg:regular_BnC} and \ref{alg:lb_BnC}. To examine the impact of inequalities (\ref{eq:c8}) on the MILP formulation of the master problem, we solve only the minimisation of total tardiness with a time limit of 900 seconds twice with and without (\ref{eq:c8}), as that objective shows weaker lower bounds.

The results regarding Algorithms \ref{alg:regular_BnC} and \ref{alg:lb_BnC} are presented in Table \ref{tab:results}. Column $R$ indicates the number of renewable resources, `Alg' indicates Algorithm (\ref{alg:regular_BnC} or \ref{alg:lb_BnC}), and columns $\sum_{j} C_{j}$ and $\sum_{j}T_{j}$ denote the total completion times and the total tardiness, respectively. `Time' indicates the elapsed time (3600 seconds being the time limit), `N' is the number of integer solutions of the master problem found during `Time'. `LB' and `UB' are the lower and upper bound respectively, while the optimality `Gap' equals $\frac{\text{UB} - \text{LB}}{\text{UB}} \%$. The values in bold denote the best-of-two value for `LB', `UB' and `Gap' between Algorithms \ref{alg:regular_BnC} and \ref{alg:lb_BnC}, e.g., for instance $J50\_M5\_\tau 0.5\_\alpha 0$, $R = 2$ and total tardiness. \textit{876} is a tightest lower bound than \textit{869}, while \textit{972} is a better upper bound than \textit{986}).

        \begin{longtable}[c]{|c|c|c|ccccc|ccccc|}
        \hline \scriptsize
        \multirow{2}{*}{Instance}                      & \multirow{2}{*}{R} & \multirow{2}{*}{Alg} & \multicolumn{5}{c|}{$\sum_{j} C_{j}$}                                                                               & \multicolumn{5}{c|}{$\sum_{j}T_{j}$}                                                                                       \\ \cline{4-13} 
                                                       &                    &                            & \multicolumn{1}{c|}{\textbf{Time}} & \multicolumn{1}{c|}{\textbf{N}}   & \multicolumn{1}{c|}{\textbf{LB}}   & \multicolumn{1}{c|}{\textbf{UB}}    & \textbf{Gap}   & \multicolumn{1}{c|}{\textbf{Time}} & \multicolumn{1}{c|}{\textbf{N}}   & \multicolumn{1}{c|}{\textbf{LB}}   & \multicolumn{1}{c|}{\textbf{UB}}   & \textbf{Gap}   \\ \hline \hline
        \endfirsthead
        \endhead
        \multirow{4}{*}{$J50\_M5\_\tau0.5\_\alpha0$}   & \multirow{2}{*}{2} & 1                          & \multicolumn{1}{c|}{232}  & \multicolumn{1}{c|}{237} & \multicolumn{1}{c|}{\textbf{4058}} & \multicolumn{1}{c|}{\textbf{4058}}  & \textbf{0.00}  & \multicolumn{1}{c|}{3600} & \multicolumn{1}{c|}{93}  & \multicolumn{1}{c|}{\textbf{876}}  & \multicolumn{1}{c|}{986}  & 11.16 \\ \cline{3-13} 
                                                       &                    & 2                          & \multicolumn{1}{c|}{189}  & \multicolumn{1}{c|}{25}  & \multicolumn{1}{c|}{\textbf{4058}} & \multicolumn{1}{c|}{\textbf{4058}}  & \textbf{0.00}  & \multicolumn{1}{c|}{3600} & \multicolumn{1}{c|}{40}  & \multicolumn{1}{c|}{869}  & \multicolumn{1}{c|}{\textbf{972}}  & \textbf{10.60} \\ \cline{2-13} 
                                                       & \multirow{2}{*}{3} & 1                          & \multicolumn{1}{c|}{60}   & \multicolumn{1}{c|}{22}  & \multicolumn{1}{c|}{\textbf{4051}} & \multicolumn{1}{c|}{\textbf{4051}}  & \textbf{0.00}  & \multicolumn{1}{c|}{3600} & \multicolumn{1}{c|}{22}  & \multicolumn{1}{c|}{873}  & \multicolumn{1}{c|}{975}  & 10.46 \\ \cline{3-13} 
                                                       &                    & 2                          & \multicolumn{1}{c|}{59}   & \multicolumn{1}{c|}{21}  & \multicolumn{1}{c|}{\textbf{4051}} & \multicolumn{1}{c|}{\textbf{4051}}  & \textbf{0.00}  & \multicolumn{1}{c|}{3600} & \multicolumn{1}{c|}{29}  & \multicolumn{1}{c|}{\textbf{879}}  & \multicolumn{1}{c|}{\textbf{960}}  & \textbf{8.44}  \\ \hline \hline
        \multirow{4}{*}{$J50\_M5\_\tau0.5\_\alpha1$}   & \multirow{2}{*}{2} & 1                          & \multicolumn{1}{c|}{3600} & \multicolumn{1}{c|}{222} & \multicolumn{1}{c|}{5589} & \multicolumn{1}{c|}{6249}  & 10.56 & \multicolumn{1}{c|}{3600} & \multicolumn{1}{c|}{209} & \multicolumn{1}{c|}{1105} & \multicolumn{1}{c|}{1806} & 38.82 \\ \cline{3-13} 
                                                       &                    & 2                          & \multicolumn{1}{c|}{1536} & \multicolumn{1}{c|}{39}  & \multicolumn{1}{c|}{\textbf{6153}} & \multicolumn{1}{c|}{\textbf{6153}}  & \textbf{0.00}  & \multicolumn{1}{c|}{3600} & \multicolumn{1}{c|}{51}  & \multicolumn{1}{c|}{\textbf{1116}} & \multicolumn{1}{c|}{\textbf{1690}} & \textbf{33.96} \\ \cline{2-13} 
                                                       & \multirow{2}{*}{3} & 1                          & \multicolumn{1}{c|}{3600} & \multicolumn{1}{c|}{288} & \multicolumn{1}{c|}{5792} & \multicolumn{1}{c|}{\textbf{5887}}  & 1.61  & \multicolumn{1}{c|}{3600} & \multicolumn{1}{c|}{277} & \multicolumn{1}{c|}{1110} & \multicolumn{1}{c|}{1482} & 25.10 \\ \cline{3-13} 
                                                       &                    & 2                          & \multicolumn{1}{c|}{821}  & \multicolumn{1}{c|}{33}  & \multicolumn{1}{c|}{\textbf{5887}} & \multicolumn{1}{c|}{\textbf{5887}}  & \textbf{0.00}  & \multicolumn{1}{c|}{3600} & \multicolumn{1}{c|}{44}  & \multicolumn{1}{c|}{\textbf{1116}} & \multicolumn{1}{c|}{\textbf{1471}} & \textbf{24.13} \\ \hline \hline
        \multirow{4}{*}{$J50\_M5\_\tau0.5\_\alpha2$}   & \multirow{2}{*}{2} & 1                          & \multicolumn{1}{c|}{3600} & \multicolumn{1}{c|}{355} & \multicolumn{1}{c|}{5925} & \multicolumn{1}{c|}{6386}  & 7.22  & \multicolumn{1}{c|}{3600} & \multicolumn{1}{c|}{219} & \multicolumn{1}{c|}{1346} & \multicolumn{1}{c|}{2199} & 38.79 \\ \cline{3-13} 
                                                       &                    & 2                          & \multicolumn{1}{c|}{3600} & \multicolumn{1}{c|}{63}  & \multicolumn{1}{c|}{\textbf{5937}} & \multicolumn{1}{c|}{\textbf{6295}}  & \textbf{5.69}  & \multicolumn{1}{c|}{3600} & \multicolumn{1}{c|}{54}  & \multicolumn{1}{c|}{\textbf{1360}} & \multicolumn{1}{c|}{\textbf{1835}} & \textbf{25.89} \\ \cline{2-13} 
                                                       & \multirow{2}{*}{3} & 1                          & \multicolumn{1}{c|}{3600} & \multicolumn{1}{c|}{248} & \multicolumn{1}{c|}{5933} & \multicolumn{1}{c|}{6163}  & 3.73  & \multicolumn{1}{c|}{3600} & \multicolumn{1}{c|}{115} & \multicolumn{1}{c|}{1347} & \multicolumn{1}{c|}{1784} & 24.50 \\ \cline{3-13} 
                                                       &                    & 2                          & \multicolumn{1}{c|}{3600} & \multicolumn{1}{c|}{55}  & \multicolumn{1}{c|}{\textbf{5955}} & \multicolumn{1}{c|}{\textbf{6119}}  & \textbf{2.68}  & \multicolumn{1}{c|}{3600} & \multicolumn{1}{c|}{51}  & \multicolumn{1}{c|}{\textbf{1368}} & \multicolumn{1}{c|}{\textbf{1702}} & \textbf{19.62} \\ \hline \hline
        \multirow{4}{*}{$J50\_M10\_\tau0.5\_\alpha0$}  & \multirow{2}{*}{4} & 1                          & \multicolumn{1}{c|}{115}  & \multicolumn{1}{c|}{185} & \multicolumn{1}{c|}{\textbf{1864}} & \multicolumn{1}{c|}{\textbf{1864}}  & \textbf{0.00}  & \multicolumn{1}{c|}{3600} & \multicolumn{1}{c|}{18}  & \multicolumn{1}{c|}{\textbf{558}}  & \multicolumn{1}{c|}{\textbf{562}}  & \textbf{0.71}  \\ \cline{3-13} 
                                                       &                    & 2                          & \multicolumn{1}{c|}{124}  & \multicolumn{1}{c|}{15}  & \multicolumn{1}{c|}{\textbf{1864}} & \multicolumn{1}{c|}{\textbf{1864}}  & \textbf{0.00}  & \multicolumn{1}{c|}{3600} & \multicolumn{1}{c|}{13}  & \multicolumn{1}{c|}{557}  & \multicolumn{1}{c|}{\textbf{562}}  & 0.89  \\ \cline{2-13} 
                                                       & \multirow{2}{*}{6} & 1                          & \multicolumn{1}{c|}{38}   & \multicolumn{1}{c|}{10}  & \multicolumn{1}{c|}{\textbf{1858}} & \multicolumn{1}{c|}{\textbf{1858}}  & \textbf{0.00}  & \multicolumn{1}{c|}{3600} & \multicolumn{1}{c|}{13}  & \multicolumn{1}{c|}{\textbf{560}}  & \multicolumn{1}{c|}{\textbf{562}}  & \textbf{0.36}  \\ \cline{3-13} 
                                                       &                    & 2                          & \multicolumn{1}{c|}{37}   & \multicolumn{1}{c|}{14}  & \multicolumn{1}{c|}{\textbf{1858}} & \multicolumn{1}{c|}{\textbf{1858}}  & \textbf{0.00}  & \multicolumn{1}{c|}{3600} & \multicolumn{1}{c|}{15}  & \multicolumn{1}{c|}{557}  & \multicolumn{1}{c|}{\textbf{562}}  & 0.89  \\ \hline \hline
        \multirow{4}{*}{$J50\_M10\_\tau0.5\_\alpha1$}  & \multirow{2}{*}{4} & 1                          & \multicolumn{1}{c|}{3600} & \multicolumn{1}{c|}{225} & \multicolumn{1}{c|}{2536} & \multicolumn{1}{c|}{2884}  & 12.07 & \multicolumn{1}{c|}{3600} & \multicolumn{1}{c|}{179} & \multicolumn{1}{c|}{694}  & \multicolumn{1}{c|}{964}  & 28.01 \\ \cline{3-13} 
                                                       &                    & 2                          & \multicolumn{1}{c|}{1108} & \multicolumn{1}{c|}{35}  & \multicolumn{1}{c|}{\textbf{2693}} & \multicolumn{1}{c|}{\textbf{2832}}  & \textbf{4.91}  & \multicolumn{1}{c|}{3600} & \multicolumn{1}{c|}{26}  & \multicolumn{1}{c|}{\textbf{706}}  & \multicolumn{1}{c|}{\textbf{946}}  & \textbf{25.37} \\ \cline{2-13} 
                                                       & \multirow{2}{*}{6} & 1                          & \multicolumn{1}{c|}{3600} & \multicolumn{1}{c|}{132} & \multicolumn{1}{c|}{2670} & \multicolumn{1}{c|}{\textbf{2700}}  & 1.11  & \multicolumn{1}{c|}{3600} & \multicolumn{1}{c|}{77}  & \multicolumn{1}{c|}{\textbf{704}}  & \multicolumn{1}{c|}{\textbf{828}}  & \textbf{14.98} \\ \cline{3-13} 
                                                       &                    & 2                          & \multicolumn{1}{c|}{463}  & \multicolumn{1}{c|}{21}  & \multicolumn{1}{c|}{\textbf{2700}} & \multicolumn{1}{c|}{\textbf{2700}}  & \textbf{0.00}  & \multicolumn{1}{c|}{3600} & \multicolumn{1}{c|}{29}  & \multicolumn{1}{c|}{699}  & \multicolumn{1}{c|}{838}  & 16.59 \\ \hline \hline
        \multirow{4}{*}{$J50\_M10\_\tau0.5\_\alpha2$}  & \multirow{2}{*}{4} & 1                          & \multicolumn{1}{c|}{3600} & \multicolumn{1}{c|}{244} & \multicolumn{1}{c|}{2700} & \multicolumn{1}{c|}{2846}  & 5.13  & \multicolumn{1}{c|}{3600} & \multicolumn{1}{c|}{168} & \multicolumn{1}{c|}{814}  & \multicolumn{1}{c|}{1119} & 27.26 \\ \cline{3-13} 
                                                       &                    & 2                          & \multicolumn{1}{c|}{750}  & \multicolumn{1}{c|}{21}  & \multicolumn{1}{c|}{\textbf{2829}} & \multicolumn{1}{c|}{\textbf{2829}}  & \textbf{0.00}  & \multicolumn{1}{c|}{3600} & \multicolumn{1}{c|}{36}  & \multicolumn{1}{c|}{\textbf{818}}  & \multicolumn{1}{c|}{\textbf{1010}} & \textbf{19.01} \\ \cline{2-13} 
                                                       & \multirow{2}{*}{6} & 1                          & \multicolumn{1}{c|}{3600} & \multicolumn{1}{c|}{47}  & \multicolumn{1}{c|}{2755} & \multicolumn{1}{c|}{\textbf{2758}}  & 0.01  & \multicolumn{1}{c|}{3600} & \multicolumn{1}{c|}{21}  & \multicolumn{1}{c|}{816}  & \multicolumn{1}{c|}{\textbf{992}}  & 17.74 \\ \cline{3-13} 
                                                       &                    & 2                          & \multicolumn{1}{c|}{1620} & \multicolumn{1}{c|}{28}  & \multicolumn{1}{c|}{\textbf{2758}} & \multicolumn{1}{c|}{\textbf{2758}}  & \textbf{0.00}  & \multicolumn{1}{c|}{3600} & \multicolumn{1}{c|}{20}  & \multicolumn{1}{c|}{\textbf{818}}  & \multicolumn{1}{c|}{993}  & \textbf{17.62} \\ \hline \hline
        \multirow{4}{*}{$J50\_M5\_\tau0.8\_\alpha0$}   & \multirow{2}{*}{2} & 1                          & \multicolumn{1}{c|}{112}  & \multicolumn{1}{c|}{59}  & \multicolumn{1}{c|}{\textbf{4752}} & \multicolumn{1}{c|}{\textbf{4752}}  & \textbf{0.00}  & \multicolumn{1}{c|}{3600} & \multicolumn{1}{c|}{194} & \multicolumn{1}{c|}{2817} & \multicolumn{1}{c|}{2873} & 1.95  \\ \cline{3-13} 
                                                       &                    & 2                          & \multicolumn{1}{c|}{141}  & \multicolumn{1}{c|}{13}  & \multicolumn{1}{c|}{\textbf{4752}} & \multicolumn{1}{c|}{\textbf{4752}}  & \textbf{0.00}  & \multicolumn{1}{c|}{3600} & \multicolumn{1}{c|}{29}  & \multicolumn{1}{c|}{\textbf{2826}} & \multicolumn{1}{c|}{\textbf{2866}} & \textbf{1.40}  \\ \cline{2-13} 
                                                       & \multirow{2}{*}{3} & 1                          & \multicolumn{1}{c|}{75}   & \multicolumn{1}{c|}{15}  & \multicolumn{1}{c|}{\textbf{4746}} & \multicolumn{1}{c|}{\textbf{4746}}  & \textbf{0.00}  & \multicolumn{1}{c|}{3600} & \multicolumn{1}{c|}{24}  & \multicolumn{1}{c|}{\textbf{2829}} & \multicolumn{1}{c|}{\textbf{2852}} & \textbf{0.81}  \\ \cline{3-13} 
                                                       &                    & 2                          & \multicolumn{1}{c|}{69}   & \multicolumn{1}{c|}{11}  & \multicolumn{1}{c|}{\textbf{4746}} & \multicolumn{1}{c|}{\textbf{4746}}  & \textbf{0.00}  & \multicolumn{1}{c|}{3600} & \multicolumn{1}{c|}{18}  & \multicolumn{1}{c|}{2822} & \multicolumn{1}{c|}{2851} & 1.02  \\ \hline \hline
        \multirow{4}{*}{$J50\_M5\_\tau0.8\_\alpha1$}   & \multirow{2}{*}{2} & 1                          & \multicolumn{1}{c|}{3600} & \multicolumn{1}{c|}{219} & \multicolumn{1}{c|}{5411} & \multicolumn{1}{c|}{6412}  & 15.61 & \multicolumn{1}{c|}{3600} & \multicolumn{1}{c|}{203} & \multicolumn{1}{c|}{3341} & \multicolumn{1}{c|}{4154} & 19.57 \\ \cline{3-13} 
                                                       &                    & 2                          & \multicolumn{1}{c|}{3600} & \multicolumn{1}{c|}{48}  & \multicolumn{1}{c|}{\textbf{5817}} & \multicolumn{1}{c|}{\textbf{6234}}  & \textbf{6.69}  & \multicolumn{1}{c|}{3600} & \multicolumn{1}{c|}{47}  & \multicolumn{1}{c|}{\textbf{3431}} & \multicolumn{1}{c|}{\textbf{4078}} & \textbf{15.87} \\ \cline{2-13} 
                                                       & \multirow{2}{*}{3} & 1                          & \multicolumn{1}{c|}{3600} & \multicolumn{1}{c|}{526} & \multicolumn{1}{c|}{\textbf{5807}} & \multicolumn{1}{c|}{\textbf{5929}}  & \textbf{2.06}  & \multicolumn{1}{c|}{3600} & \multicolumn{1}{c|}{234} & \multicolumn{1}{c|}{3418} & \multicolumn{1}{c|}{\textbf{3713}} & 7.95  \\ \cline{3-13} 
                                                       &                    & 2                          & \multicolumn{1}{c|}{3600} & \multicolumn{1}{c|}{33}  & \multicolumn{1}{c|}{5659} & \multicolumn{1}{c|}{5956}  & 4.99  & \multicolumn{1}{c|}{3600} & \multicolumn{1}{c|}{27}  & \multicolumn{1}{c|}{\textbf{3439}} & \multicolumn{1}{c|}{3717} & \textbf{7.48}  \\ \hline \hline
        \multirow{4}{*}{$J50\_M5\_\tau0.8\_\alpha2$}   & \multirow{2}{*}{2} & 1                          & \multicolumn{1}{c|}{3600} & \multicolumn{1}{c|}{249} & \multicolumn{1}{c|}{\textbf{5661}} & \multicolumn{1}{c|}{6289}  & 9.99 & \multicolumn{1}{c|}{3600} & \multicolumn{1}{c|}{210} & \multicolumn{1}{c|}{3476} & \multicolumn{1}{c|}{4343} & 19.96 \\ \cline{3-13} 
                                                       &                    & 2                          & \multicolumn{1}{c|}{3600} & \multicolumn{1}{c|}{40}  & \multicolumn{1}{c|}{5637} & \multicolumn{1}{c|}{\textbf{6088}}  & \textbf{7.41}  & \multicolumn{1}{c|}{3600} & \multicolumn{1}{c|}{51}  & \multicolumn{1}{c|}{\textbf{3519}} & \multicolumn{1}{c|}{\textbf{3968}} & \textbf{11.32} \\ \cline{2-13} 
                                                       & \multirow{2}{*}{3} & 1                          & \multicolumn{1}{c|}{3600} & \multicolumn{1}{c|}{303} & \multicolumn{1}{c|}{5661} & \multicolumn{1}{c|}{5912}  & 4.25  & \multicolumn{1}{c|}{3600} & \multicolumn{1}{c|}{76}  & \multicolumn{1}{c|}{\textbf{3544}} & \multicolumn{1}{c|}{3809} & \textbf{6.96}  \\ \cline{3-13} 
                                                       &                    & 2                          & \multicolumn{1}{c|}{3600} & \multicolumn{1}{c|}{29}  & \multicolumn{1}{c|}{\textbf{5725}} & \multicolumn{1}{c|}{\textbf{5882}}  & \textbf{2.67}  & \multicolumn{1}{c|}{3600} & \multicolumn{1}{c|}{52}  & \multicolumn{1}{c|}{3526} & \multicolumn{1}{c|}{\textbf{3801}} & 7.23  \\ \hline \hline
        \multirow{4}{*}{$J50\_M10\_\tau0.8\_\alpha0$}  & \multirow{2}{*}{4} & 1                          & \multicolumn{1}{c|}{79}   & \multicolumn{1}{c|}{23}  & \multicolumn{1}{c|}{\textbf{2393}} & \multicolumn{1}{c|}{\textbf{2393}}  & \textbf{0.00}  & \multicolumn{1}{c|}{386}  & \multicolumn{1}{c|}{27}  & \multicolumn{1}{c|}{\textbf{1622}} & \multicolumn{1}{c|}{\textbf{1622}} & \textbf{0.00}  \\ \cline{3-13} 
                                                       &                    & 2                          & \multicolumn{1}{c|}{108}  & \multicolumn{1}{c|}{13}  & \multicolumn{1}{c|}{\textbf{2393}} & \multicolumn{1}{c|}{\textbf{2393}}  & \textbf{0.00}  & \multicolumn{1}{c|}{375}  & \multicolumn{1}{c|}{14}  & \multicolumn{1}{c|}{\textbf{1622}} & \multicolumn{1}{c|}{\textbf{1622}} & \textbf{0.00}  \\ \cline{2-13} 
                                                       & \multirow{2}{*}{6} & 1                          & \multicolumn{1}{c|}{66}   & \multicolumn{1}{c|}{15}  & \multicolumn{1}{c|}{\textbf{2391}} & \multicolumn{1}{c|}{\textbf{2391}}  & \textbf{0.00}  & \multicolumn{1}{c|}{264}  & \multicolumn{1}{c|}{6}   & \multicolumn{1}{c|}{\textbf{1619}} & \multicolumn{1}{c|}{\textbf{1619}} & \textbf{0.00}  \\ \cline{3-13} 
                                                       &                    & 2                          & \multicolumn{1}{c|}{60}   & \multicolumn{1}{c|}{12}  & \multicolumn{1}{c|}{\textbf{2391}} & \multicolumn{1}{c|}{\textbf{2391}}  & \textbf{0.00}  & \multicolumn{1}{c|}{262}  & \multicolumn{1}{c|}{6}   & \multicolumn{1}{c|}{\textbf{1619}} & \multicolumn{1}{c|}{\textbf{1619}} & \textbf{0.00}  \\ \hline \hline
        \multirow{4}{*}{$J50\_M10\_\tau0.8\_\alpha1$}  & \multirow{2}{*}{4} & 1                          & \multicolumn{1}{c|}{3600} & \multicolumn{1}{c|}{226} & \multicolumn{1}{c|}{2684} & \multicolumn{1}{c|}{2851}  & 5.86  & \multicolumn{1}{c|}{3600} & \multicolumn{1}{c|}{212} & \multicolumn{1}{c|}{1711} & \multicolumn{1}{c|}{1900} & 9.95  \\ \cline{3-13} 
                                                       &                    & 2                          & \multicolumn{1}{c|}{806}  & \multicolumn{1}{c|}{29}  & \multicolumn{1}{c|}{\textbf{2833}} & \multicolumn{1}{c|}{\textbf{2833}}  & \textbf{0.00}  & \multicolumn{1}{c|}{670}  & \multicolumn{1}{c|}{17}  & \multicolumn{1}{c|}{\textbf{1859}} & \multicolumn{1}{c|}{\textbf{1859}} & \textbf{0.00}  \\ \cline{2-13} 
                                                       & \multirow{2}{*}{6} & 1                          & \multicolumn{1}{c|}{195}  & \multicolumn{1}{c|}{103} & \multicolumn{1}{c|}{\textbf{2722}} & \multicolumn{1}{c|}{\textbf{2722}}  & \textbf{0.00}  & \multicolumn{1}{c|}{1138} & \multicolumn{1}{c|}{65}  & \multicolumn{1}{c|}{\textbf{1753}} & \multicolumn{1}{c|}{\textbf{1753}} & \textbf{0.00}  \\ \cline{3-13} 
                                                       &                    & 2                          & \multicolumn{1}{c|}{394}  & \multicolumn{1}{c|}{25}  & \multicolumn{1}{c|}{\textbf{2722}} & \multicolumn{1}{c|}{\textbf{2722}}  & \textbf{0.00}  & \multicolumn{1}{c|}{362}  & \multicolumn{1}{c|}{23}  & \multicolumn{1}{c|}{\textbf{1753}} & \multicolumn{1}{c|}{\textbf{1753}} & \textbf{0.00}  \\ \hline \hline
        \multirow{4}{*}{$J50\_M10\_\tau0.8\_\alpha2$}  & \multirow{2}{*}{4} & 1                          & \multicolumn{1}{c|}{3600} & \multicolumn{1}{c|}{215} & \multicolumn{1}{c|}{2749} & \multicolumn{1}{c|}{2990}  & 8.06  & \multicolumn{1}{c|}{3600} & \multicolumn{1}{c|}{182} & \multicolumn{1}{c|}{1860} & \multicolumn{1}{c|}{2209} & 15.80 \\ \cline{3-13} 
                                                       &                    & 2                          & \multicolumn{1}{c|}{3600} & \multicolumn{1}{c|}{46}  & \multicolumn{1}{c|}{\textbf{2791}} & \multicolumn{1}{c|}{\textbf{2935}}  & \textbf{4.91}  & \multicolumn{1}{c|}{3600} & \multicolumn{1}{c|}{45}  & \multicolumn{1}{c|}{\textbf{1900}} & \multicolumn{1}{c|}{\textbf{2137}} & \textbf{11.09} \\ \cline{2-13} 
                                                       & \multirow{2}{*}{6} & 1                          & \multicolumn{1}{c|}{3600} & \multicolumn{1}{c|}{93}  & \multicolumn{1}{c|}{\textbf{2781}} & \multicolumn{1}{c|}{\textbf{2834}}  & \textbf{1.87}  & \multicolumn{1}{c|}{3600} & \multicolumn{1}{c|}{60}  & \multicolumn{1}{c|}{\textbf{1906}} & \multicolumn{1}{c|}{2018} & \textbf{5.55}  \\ \cline{3-13} 
                                                       &                    & 2                          & \multicolumn{1}{c|}{3600} & \multicolumn{1}{c|}{29}  & \multicolumn{1}{c|}{2769} & \multicolumn{1}{c|}{2849}  & 2.81  & \multicolumn{1}{c|}{3600} & \multicolumn{1}{c|}{37}  & \multicolumn{1}{c|}{1900} & \multicolumn{1}{c|}{\textbf{2017}} & 5.80  \\ \hline \hline
        \multirow{4}{*}{$J100\_M10\_\tau0.5\_\alpha0$} & \multirow{2}{*}{4} & 1                          & \multicolumn{1}{c|}{3600} & \multicolumn{1}{c|}{46}  & \multicolumn{1}{c|}{\textbf{7284}} & \multicolumn{1}{c|}{7362}  & \textbf{1.06}  & \multicolumn{1}{c|}{3600} & \multicolumn{1}{c|}{23}  & \multicolumn{1}{c|}{\textbf{1545}} & \multicolumn{1}{c|}{1826} & \textbf{15.39} \\ \cline{3-13} 
                                                       &                    & 2                          & \multicolumn{1}{c|}{3600} & \multicolumn{1}{c|}{32}  & \multicolumn{1}{c|}{7274} & \multicolumn{1}{c|}{\textbf{7352}}  & \textbf{1.06}  & \multicolumn{1}{c|}{3600} & \multicolumn{1}{c|}{32}  & \multicolumn{1}{c|}{1521} & \multicolumn{1}{c|}{\textbf{1814}} & 16.15 \\ \cline{2-13} 
                                                       & \multirow{2}{*}{6} & 1                          & \multicolumn{1}{c|}{3600} & \multicolumn{1}{c|}{32}  & \multicolumn{1}{c|}{\textbf{7250}} & \multicolumn{1}{c|}{\textbf{7339}}  & \textbf{1.21}  & \multicolumn{1}{c|}{3600} & \multicolumn{1}{c|}{18}  & \multicolumn{1}{c|}{\textbf{1553}} & \multicolumn{1}{c|}{\textbf{1806}} & \textbf{14.01} \\ \cline{3-13} 
                                                       &                    & 2                          & \multicolumn{1}{c|}{3600} & \multicolumn{1}{c|}{19}  & \multicolumn{1}{c|}{7236} & \multicolumn{1}{c|}{\textbf{7339}}  & 1.40  & \multicolumn{1}{c|}{3600} & \multicolumn{1}{c|}{17}  & \multicolumn{1}{c|}{1551} & \multicolumn{1}{c|}{1825} & 15.01 \\ \hline \hline
        \multirow{4}{*}{$J100\_M10\_\tau0.5\_\alpha1$} & \multirow{2}{*}{4} & 1                          & \multicolumn{1}{c|}{3600} & \multicolumn{1}{c|}{137} & \multicolumn{1}{c|}{8327} & \multicolumn{1}{c|}{10593} & 21.39 & \multicolumn{1}{c|}{3600} & \multicolumn{1}{c|}{55}  & \multicolumn{1}{c|}{1284} & \multicolumn{1}{c|}{\textbf{2863}} & \textbf{55.15} \\ \cline{3-13} 
                                                       &                    & 2                          & \multicolumn{1}{c|}{3600} & \multicolumn{1}{c|}{30}  & \multicolumn{1}{c|}{\textbf{9043}} & \multicolumn{1}{c|}{\textbf{10334}} & \textbf{12.49} & \multicolumn{1}{c|}{3600} & \multicolumn{1}{c|}{24}  & \multicolumn{1}{c|}{\textbf{1289}} & \multicolumn{1}{c|}{3030} & 57.46 \\ \cline{2-13} 
                                                       & \multirow{2}{*}{6} & 1                          & \multicolumn{1}{c|}{3600} & \multicolumn{1}{c|}{67}  & \multicolumn{1}{c|}{9210} & \multicolumn{1}{c|}{9843}  & 6.43  & \multicolumn{1}{c|}{3600} & \multicolumn{1}{c|}{32}  & \multicolumn{1}{c|}{\textbf{1298}} & \multicolumn{1}{c|}{2570} & \textbf{49.49} \\ \cline{3-13} 
                                                       &                    & 2                          & \multicolumn{1}{c|}{3600} & \multicolumn{1}{c|}{38}  & \multicolumn{1}{c|}{\textbf{9231}} & \multicolumn{1}{c|}{\textbf{9815}}  & \textbf{5.95}  & \multicolumn{1}{c|}{3600} & \multicolumn{1}{c|}{27}  & \multicolumn{1}{c|}{1279} & \multicolumn{1}{c|}{\textbf{2554}} & 49.92 \\ \hline \hline
        \multirow{4}{*}{$J100\_M10\_\tau0.5\_\alpha2$} & \multirow{2}{*}{4} & 1                          & \multicolumn{1}{c|}{3600} & \multicolumn{1}{c|}{80}  & \multicolumn{1}{c|}{9311} & \multicolumn{1}{c|}{10691} & 12.91 & \multicolumn{1}{c|}{3600} & \multicolumn{1}{c|}{50}  & \multicolumn{1}{c|}{\textbf{2023}} & \multicolumn{1}{c|}{3605} & 43.88 \\ \cline{3-13} 
                                                       &                    & 2                          & \multicolumn{1}{c|}{3600} & \multicolumn{1}{c|}{31}  & \multicolumn{1}{c|}{\textbf{9445}} & \multicolumn{1}{c|}{\textbf{10458}} & \textbf{9.69}  & \multicolumn{1}{c|}{3600} & \multicolumn{1}{c|}{29}  & \multicolumn{1}{c|}{2022} & \multicolumn{1}{c|}{\textbf{3530}} & \textbf{42.72} \\ \cline{2-13} 
                                                       & \multirow{2}{*}{6} & 1                          & \multicolumn{1}{c|}{3600} & \multicolumn{1}{c|}{65}  & \multicolumn{1}{c|}{\textbf{9507}} & \multicolumn{1}{c|}{10168} & \textbf{6.50}  & \multicolumn{1}{c|}{3600} & \multicolumn{1}{c|}{22}  & \multicolumn{1}{c|}{2029} & \multicolumn{1}{c|}{\textbf{3110}} & 34.76 \\ \cline{3-13} 
                                                       &                    & 2                          & \multicolumn{1}{c|}{3600} & \multicolumn{1}{c|}{47}  & \multicolumn{1}{c|}{9349} & \multicolumn{1}{c|}{\textbf{10063}} & 7.10  & \multicolumn{1}{c|}{3600} & \multicolumn{1}{c|}{27}  & \multicolumn{1}{c|}{\textbf{2063}} & \multicolumn{1}{c|}{3124} & \textbf{33.96} \\ \hline \hline
        \multirow{4}{*}{$J100\_M10\_\tau0.8\_\alpha0$} & \multirow{2}{*}{4} & 1                          & \multicolumn{1}{c|}{3600} & \multicolumn{1}{c|}{48}  & \multicolumn{1}{c|}{\textbf{6961}} & \multicolumn{1}{c|}{\textbf{7008}}  & \textbf{0.67}  & \multicolumn{1}{c|}{3600} & \multicolumn{1}{c|}{23}  & \multicolumn{1}{c|}{4187} & \multicolumn{1}{c|}{4336} & 3.44  \\ \cline{3-13} 
                                                       &                    & 2                          & \multicolumn{1}{c|}{3600} & \multicolumn{1}{c|}{37}  & \multicolumn{1}{c|}{6952} & \multicolumn{1}{c|}{7021}  & 0.98  & \multicolumn{1}{c|}{3600} & \multicolumn{1}{c|}{27}  & \multicolumn{1}{c|}{\textbf{4194}} & \multicolumn{1}{c|}{\textbf{4322}} & \textbf{2.96}  \\ \cline{2-13} 
                                                       & \multirow{2}{*}{6} & 1                          & \multicolumn{1}{c|}{3600} & \multicolumn{1}{c|}{28}  & \multicolumn{1}{c|}{6943} & \multicolumn{1}{c|}{\textbf{7002}}  & 0.84  & \multicolumn{1}{c|}{3600} & \multicolumn{1}{c|}{13}  & \multicolumn{1}{c|}{\textbf{4156}} & \multicolumn{1}{c|}{\textbf{4333}} & \textbf{4.08}  \\ \cline{3-13} 
                                                       &                    & 2                          & \multicolumn{1}{c|}{3600} & \multicolumn{1}{c|}{26}  & \multicolumn{1}{c|}{\textbf{6953}} & \multicolumn{1}{c|}{7009}  & \textbf{0.80}  & \multicolumn{1}{c|}{3600} & \multicolumn{1}{c|}{12}  & \multicolumn{1}{c|}{4151} & \multicolumn{1}{c|}{4336} & 4.27  \\ \hline \hline
        \multirow{4}{*}{$J100\_M10\_\tau0.8\_\alpha1$} & \multirow{2}{*}{4} & 1                          & \multicolumn{1}{c|}{3600} & \multicolumn{1}{c|}{93}  & \multicolumn{1}{c|}{8155} & \multicolumn{1}{c|}{\textbf{9694}}  & 15.88 & \multicolumn{1}{c|}{3600} & \multicolumn{1}{c|}{49}  & \multicolumn{1}{c|}{\textbf{4379}} & \multicolumn{1}{c|}{6283} & 30.30 \\ \cline{3-13} 
                                                       &                    & 2                          & \multicolumn{1}{c|}{3600} & \multicolumn{1}{c|}{25}  & \multicolumn{1}{c|}{\textbf{8740}} & \multicolumn{1}{c|}{9727}  & \textbf{10.15} & \multicolumn{1}{c|}{3600} & \multicolumn{1}{c|}{15}  & \multicolumn{1}{c|}{\textbf{4379}} & \multicolumn{1}{c|}{\textbf{6238}} & \textbf{29.80} \\ \cline{2-13} 
                                                       & \multirow{2}{*}{6} & 1                          & \multicolumn{1}{c|}{3600} & \multicolumn{1}{c|}{53}  & \multicolumn{1}{c|}{8736} & \multicolumn{1}{c|}{\textbf{9337}}  & 6.44  & \multicolumn{1}{c|}{3600} & \multicolumn{1}{c|}{31}  & \multicolumn{1}{c|}{\textbf{4411}} & \multicolumn{1}{c|}{5852} & 24.62 \\ \cline{3-13} 
                                                       &                    & 2                          & \multicolumn{1}{c|}{3600} & \multicolumn{1}{c|}{35}  & \multicolumn{1}{c|}{\textbf{8788}} & \multicolumn{1}{c|}{9386}  & \textbf{6.37}  & \multicolumn{1}{c|}{3600} & \multicolumn{1}{c|}{19}  & \multicolumn{1}{c|}{4386} & \multicolumn{1}{c|}{\textbf{5801}} & \textbf{24.39} \\ \hline \hline
        \multirow{4}{*}{$J100\_M10\_\tau0.8\_\alpha2$} & \multirow{2}{*}{4} & 1                          & \multicolumn{1}{c|}{3600} & \multicolumn{1}{c|}{84}  & \multicolumn{1}{c|}{8667} & \multicolumn{1}{c|}{10236} & 15.33 & \multicolumn{1}{c|}{3600} & \multicolumn{1}{c|}{22}  & \multicolumn{1}{c|}{\textbf{5045}} & \multicolumn{1}{c|}{7498} & 32.72 \\ \cline{3-13} 
                                                       &                    & 2                          & \multicolumn{1}{c|}{3600} & \multicolumn{1}{c|}{34}  & \multicolumn{1}{c|}{\textbf{8882}} & \multicolumn{1}{c|}{\textbf{10180}} & \textbf{12.75} & \multicolumn{1}{c|}{3600} & \multicolumn{1}{c|}{9}   & \multicolumn{1}{c|}{5024} & \multicolumn{1}{c|}{\textbf{7223}} & \textbf{30.44} \\ \cline{2-13} 
                                                       & \multirow{2}{*}{6} & 1                          & \multicolumn{1}{c|}{3600} & \multicolumn{1}{c|}{51}  & \multicolumn{1}{c|}{\textbf{8982}} & \multicolumn{1}{c|}{9962}  & 9.84  & \multicolumn{1}{c|}{3600} & \multicolumn{1}{c|}{12}  & \multicolumn{1}{c|}{5040} & \multicolumn{1}{c|}{6803} & 25.92 \\ \cline{3-13} 
                                                       &                    & 2                          & \multicolumn{1}{c|}{3600} & \multicolumn{1}{c|}{29}  & \multicolumn{1}{c|}{8976} & \multicolumn{1}{c|}{\textbf{9751}}  & \textbf{7.95}  & \multicolumn{1}{c|}{3600} & \multicolumn{1}{c|}{6}   & \multicolumn{1}{c|}{\textbf{5079}} & \multicolumn{1}{c|}{\textbf{6724}} & \textbf{24.46} \\ \hline
        \caption{Comparative results of Algorithms \ref{alg:regular_BnC} and \ref{alg:lb_BnC}}
        \label{tab:results}\\
        \end{longtable}

Regarding total completion time, we observe that the suggested LBBD improvements through local branching exhibit a significant impact on the reduction of optimality gap or for closing the gap faster, for all instances of 50 jobs. A significant number of instances is solved to optimality, several of which only after the inclusion of local branching cuts. We also notice that the reduction is larger for more difficult instances, i.e., instances of larger setup times and fewer resources like $J50\_M10\_\tau 0.5\_\alpha 1$, $R = 4$ and $J50\_M5\_\tau 0.8\_\alpha 1$, $R = 2$. The drastic decrease in the number of integer solutions encountered reflects both the speed-up accomplished and the additional effort required per solution to explore the neighbourhood.
        
For instances of $100$ jobs, the time limit is not adequate to reach optimality, with gaps rising up for $\alpha \in \{1,2\}.$ The performance of Local Branching slightly declines, not surprisingly, as the number of internal swaps, hence the size of the neighbourhoods and their exploration, increases. However, for all instances of $\alpha = 1$ and $2$ (i.e., with large setup times, compared to the group of $\alpha = 0$), Algorithm \ref{alg:lb_BnC} provides better results in terms of optimality gap. For instances of $\alpha = 0$, few improvements could be expected, as both algorithms compute solutions of $~1\%$ gap. Overall, the LBBD and Branch-and-Check we introduce works well and is further strengthened by the tighter Benders' cuts obtained through local branching.

Total tardiness is quite more challenging, with increased  optimality gaps within the 1-hour limit. For instances of $50$ jobs, Algorithm \ref{alg:lb_BnC} steadily dominates Algorithm \ref{alg:regular_BnC}, with the exception of a small number of instances (e.g., $J50\_M10\_\tau 0.5\_\alpha 1$, $R = 6$ and $J50\_M5\_\tau 0.8 \_\alpha 2$, $R = 3$). For instances of $100$ jobs, Local Branching shows improved gaps for $\alpha = 2$, while none of the algorithms shows clear dominance for $\alpha = 0$ or $1.$ In general, the differentiation noticed is less significant than for total completion times. This could be attributed on one hand to the non-linearity of total tardiness and hence the weaker bounds obtained by its linearisation: even if the UB provided by Algorithm \ref{alg:lb_BnC} decreases, the percentage of decrease is understated if LB is quite distant. More time could of course lead to smaller gaps, or more clear improvements accomplished by Algorithm \ref{alg:lb_BnC}. 

Let us now solve the MILP of $\mathcal{M}$ for the total tardiness objective with and without the inclusion of valid inequalities (\ref{eq:c8}), to show their impact on tightening the lower bound of the problem. In Table \ref{tab:milp_results}, we present the values of `LB', `UB' and `Gap' for the MILP of the master problem, before (`no inequalities (\ref{eq:c8})') and after (`with inequalities (\ref{eq:c8})') the inclusion of (\ref{eq:c8}).

        \begin{table}[H]
            \centering \scriptsize
            \resizebox{\textwidth}{!}{%
            \begin{tabular}{|c||cccc||cccc|}
            \hline
            \multirow{2}{*}{\textbf{Instance}} & \multicolumn{4}{c||}{\textbf{no inequalities (\ref{eq:c8})}}                                                                         & \multicolumn{4}{c|}{\textbf{with  inequalities (\ref{eq:c8})}}                                                                            \\ \cline{2-9} 
                                               & \multicolumn{1}{c|}{\textbf{Time}} & \multicolumn{1}{c|}{\textbf{LB}}   & \multicolumn{1}{c|}{\textbf{UB}}   & \textbf{Gap}  & \multicolumn{1}{c|}{\textbf{Time}} & \multicolumn{1}{c|}{\textbf{LB}}   & \multicolumn{1}{c|}{\textbf{UB}}   & \textbf{Gap}   \\ \hline
            $J50\_M5\_\tau0.5\_\alpha0$        & \multicolumn{1}{c|}{900}           & \multicolumn{1}{c|}{874}           & \multicolumn{1}{c|}{972}           & 10.08         & \multicolumn{1}{c|}{900}           & \multicolumn{1}{c|}{\textbf{900}}  & \multicolumn{1}{c|}{\textbf{966}}  & \textbf{6.83}  \\ \hline
            $J50\_M5\_\tau0.8\_\alpha0$        & \multicolumn{1}{c|}{900}           & \multicolumn{1}{c|}{2830}          & \multicolumn{1}{c|}{2852}          & 0.77          & \multicolumn{1}{c|}{409}           & \multicolumn{1}{c|}{\textbf{2851}} & \multicolumn{1}{c|}{\textbf{2851}} & \textbf{0.00}  \\ \hline
            $J50\_M5\_\tau0.5\_\alpha1$        & \multicolumn{1}{c|}{900}           & \multicolumn{1}{c|}{1094}          & \multicolumn{1}{c|}{1466}          & 25.38         & \multicolumn{1}{c|}{900}           & \multicolumn{1}{c|}{\textbf{1252}} & \multicolumn{1}{c|}{\textbf{1423}} & \textbf{12.02} \\ \hline
            $J50\_M5\_\tau0.8\_\alpha1$        & \multicolumn{1}{c|}{900}           & \multicolumn{1}{c|}{3384}          & \multicolumn{1}{c|}{3633}          & 6.85          & \multicolumn{1}{c|}{900}           & \multicolumn{1}{c|}{\textbf{3497}} & \multicolumn{1}{c|}{\textbf{3627}} & \textbf{3.58}  \\ \hline
            $J50\_M5\_\tau0.5\_\alpha2$        & \multicolumn{1}{c|}{900}           & \multicolumn{1}{c|}{1266}          & \multicolumn{1}{c|}{1722}          & 26.48         & \multicolumn{1}{c|}{900}           & \multicolumn{1}{c|}{\textbf{1471}} & \multicolumn{1}{c|}{\textbf{1672}} & \textbf{12.02} \\ \hline
            $J50\_M5\_\tau0.8\_\alpha2$        & \multicolumn{1}{c|}{900}           & \multicolumn{1}{c|}{3321}          & \multicolumn{1}{c|}{3797}          & 12.54         & \multicolumn{1}{c|}{900}           & \multicolumn{1}{c|}{\textbf{3567}} & \multicolumn{1}{c|}{\textbf{3736}} & \textbf{4.52}  \\ \hline
            $J50\_M10\_\tau0.5\_\alpha0$       & \multicolumn{1}{c|}{900}           & \multicolumn{1}{c|}{552}           & \multicolumn{1}{c|}{\textbf{562}}  & 1.78          & \multicolumn{1}{c|}{793}           & \multicolumn{1}{c|}{\textbf{562}}  & \multicolumn{1}{c|}{\textbf{562}}  & \textbf{0.00}  \\ \hline
            $J50\_M10\_\tau0.8\_\alpha0$       & \multicolumn{1}{c|}{410}           & \multicolumn{1}{c|}{\textbf{1619}} & \multicolumn{1}{c|}{\textbf{1619}} & \textbf{0.00} & \multicolumn{1}{c|}{228}           & \multicolumn{1}{c|}{\textbf{1619}} & \multicolumn{1}{c|}{\textbf{1619}} & \textbf{0.00}  \\ \hline
            $J50\_M10\_\tau0.5\_\alpha1$       & \multicolumn{1}{c|}{900}           & \multicolumn{1}{c|}{706}           & \multicolumn{1}{c|}{821}           & 14.01         & \multicolumn{1}{c|}{900}           & \multicolumn{1}{c|}{\textbf{743}}  & \multicolumn{1}{c|}{\textbf{813}}  & \textbf{8.61}  \\ \hline
            $J50\_M10\_\tau0.8\_\alpha1$       & \multicolumn{1}{c|}{900}           & \multicolumn{1}{c|}{1714}          & \multicolumn{1}{c|}{1748}          & 1.95          & \multicolumn{1}{c|}{253}           & \multicolumn{1}{c|}{\textbf{1747}} & \multicolumn{1}{c|}{\textbf{1747}} & \textbf{0.00}  \\ \hline
            $J50\_M10\_\tau0.5\_\alpha2$       & \multicolumn{1}{c|}{900}           & \multicolumn{1}{c|}{721}           & \multicolumn{1}{c|}{988}           & 27.02         & \multicolumn{1}{c|}{900}           & \multicolumn{1}{c|}{\textbf{849}}  & \multicolumn{1}{c|}{\textbf{961}}  & \textbf{11.65} \\ \hline
            $J50\_M10\_\tau0.8\_\alpha2$       & \multicolumn{1}{c|}{900}           & \multicolumn{1}{c|}{1757}          & \multicolumn{1}{c|}{2017}          & 12.89         & \multicolumn{1}{c|}{900}           & \multicolumn{1}{c|}{\textbf{1938}} & \multicolumn{1}{c|}{\textbf{2012}} & \textbf{3.68}  \\ \hline
            $J100\_M10\_\tau0.5\_\alpha0$      & \multicolumn{1}{c|}{900}           & \multicolumn{1}{c|}{1202}          & \multicolumn{1}{c|}{1854}          & 35.17         & \multicolumn{1}{c|}{900}           & \multicolumn{1}{c|}{\textbf{1445}} & \multicolumn{1}{c|}{\textbf{1838}} & \textbf{21.38} \\ \hline
            $J100\_M10\_\tau0.8\_\alpha0$      & \multicolumn{1}{c|}{900}           & \multicolumn{1}{c|}{3617}          & \multicolumn{1}{c|}{4379}          & 17.40         & \multicolumn{1}{c|}{900}           & \multicolumn{1}{c|}{\textbf{3949}} & \multicolumn{1}{c|}{\textbf{4385}} & \textbf{9.94}  \\ \hline
            $J100\_M10\_\tau0.5\_\alpha1$      & \multicolumn{1}{c|}{900}           & \multicolumn{1}{c|}{576}           & \multicolumn{1}{c|}{2675}          & 78.47         & \multicolumn{1}{c|}{900}           & \multicolumn{1}{c|}{\textbf{1292}} & \multicolumn{1}{c|}{\textbf{2519}} & \textbf{48.71} \\ \hline
            $J100\_M10\_\tau0.8\_\alpha1$      & \multicolumn{1}{c|}{900}           & \multicolumn{1}{c|}{3122}          & \multicolumn{1}{c|}{5934}          & 47.39         & \multicolumn{1}{c|}{900}           & \multicolumn{1}{c|}{\textbf{4365}} & \multicolumn{1}{c|}{\textbf{5881}} & \textbf{25.78} \\ \hline
            $J100\_M10\_\tau0.5\_\alpha2$      & \multicolumn{1}{c|}{900}           & \multicolumn{1}{c|}{816}           & \multicolumn{1}{c|}{3548}          & 77.00         & \multicolumn{1}{c|}{900}           & \multicolumn{1}{c|}{\textbf{2032}} & \multicolumn{1}{c|}{\textbf{3378}} & \textbf{39.85} \\ \hline
            $J100\_M10\_\tau0.8\_\alpha2$      & \multicolumn{1}{c|}{900}           & \multicolumn{1}{c|}{3145}          & \multicolumn{1}{c|}{6787}          & 53.66         & \multicolumn{1}{c|}{900}           & \multicolumn{1}{c|}{\textbf{5003}} & \multicolumn{1}{c|}{\textbf{6761}} & \textbf{26.00} \\ \hline
            \end{tabular}%
            }
            \caption{Impact of inequalities (\ref{eq:c8}) on MILP performance for the master problem $\mathcal{M}$}
            \label{tab:milp_results}
        \end{table}
        
It is clear that these valid inequalities impose significantly tighter lower bounds for all instances, except for one. Even for that instance, which is solved to optimality for both groups of experiments, the addition of (\ref{eq:c8}) almost halves the time to optimality. Constraints (\ref{eq:c8}) imply also improved upper bounds for all instances.

All generated instances and coding files are available at \url{https://github.com/i-avgerinos/local_branching}.

\section{Concluding remarks}
We have introduced an LBBD for minimisation of tardiness or completion times in elaborate scheduling environments and also examined local branching constraints as Benders cuts that rule out neighbourhoods rather than single solutions. Although our LBBD covers broader cases and objectives, we have focused on computationally challenging aspects like resource constraints and total tardiness. That is, the existence of unlimited resources (i.e., $R=|M|$) creates a more favourable setting for local branching, as neighbourhood exploration would be faster, hence generating stronger Benders cuts would remain beneficial for larger instances. The makespan objective, with or without resources, would also benefit more from local branching.

Our LBBD favours local search as the subproblem is quite restricted and the neighbourhood small. Hence, one could examine whether our approach remains effective for larger neighbourhoods of the same problem (e.g., the $8$-OPT) or for less-restricted subproblems within and beyond production scheduling. Another aspect, revealed by our computational work, is that a variable neighbourhood strategy could be employed, for example by increasing neighbourhood size during search or decreasing it with problem size.

To show the versatility of our approach, let us list the $6$-OPT and $8$-OPT neighbourhoods for the extended formulation $\mathcal{M}'$ that includes the model of Section \ref{sec:MILP}, plus constraints (\ref{eq:yvars}) and the vector of variables $y$.

Observe that $k$-OPT neighbourhoods exist only for $k$ even, as the total number of variables at value $1$ is constant. To define $6$-OPT and $8$-OPT, we define the \emph{internal $h$-swap} as the movement in which $h$ assigned to the same machine in the reference solution `rotate' their slots, i.e., given $(\bar{x}^{t}, \bar{y}^{t})$ and an (ordered) set of jobs $j_1,j_2,\ldots j_h$ all assigned to slots $i_1,i_2,\ldots,i_h$ of the same machine $m,$ an internal $h$-swap assigns job $j_l$ to slot $j_{l+1}$, addition done $mod(h-1)$ for $l=1,2,\ldots,h.$ Although for $h=2$ we have the traditional internal swaps, please observe that a $4$-swap is other than two $2$-swaps. Figure \ref{fig:non_swap} presents an internal $3$-swap.

\begin{proposition}
The $6$-OPT neighbourhood for problem $\mathcal{P}$ is the $4$-OPT neighbourhood plus solutions derived from the reference one through internal $3$-swaps. \label{prop:2}
\end{proposition}  
\proof
Let  $(\bar{x}^{t}, \bar{y}^{t})$ be the reference solution. As (\ref{eq:k'_constraint}) becomes
        \begin{equation}
            \sum_{m\in M}\sum_{j:\bar{y}^{t}_{jm}=1}(1 - y_{jm}) +  \sum_{m\in M}\sum_{(i,j):\bar{x}^{t}_{ijm} = 1}(1 - x_{ijm}) \leq 3,
            \label{eq:local_branching_3}
        \end{equation}
and $4$-OPT includes all solutions making the left-hand side of (\ref{eq:local_branching_3}) $2,$ we focus on solutions that make it $3$. Observe that whenever a variable $y_{jm}$ switches from $1$ to $0,$ some variable $x_{ijm}$ also does by (\ref{eq:yvars}) (but not vice-versa); also, some other variable $y_{jm'}$ switches from $0$ to $1$, hence pushing a variable $x_{i'jm'}$ to do the same. Such a solution would already make the left-hand side of (\ref{eq:local_branching_3}) equal to $2,$ hence excluding from $6$-OPT any solution that switches more than one job from one machine to some other. Notice also that switching a single such job, i.e., a starting-job insertion, cannot be complemented with a third variable switch to change $6$ variables in total. Thus, the only possibility is to examine a solution in which $y=y',$ and three variables receive a value $1$ in $x$ and $0$ in $x'$, three other variables doing the opposite. It can easily be seen that internal $3$-swaps yield exactly these solutions.
\endproof

The following can be shown through similar arguments. 
\begin{proposition}
The $8$-OPT neighbourhood for problem $\mathcal{P}$ is the $6$-OPT neighbourhood plus solutions derived from the reference one through one of the following (i) two internal swaps, (ii) two first-job insertions, (iii) one internal swap and one first-job insertion, (iv) an internal $4$-swap \label{prop:3}, (v) an external swap.
\end{proposition} 


In broader theoretical terms, it could be fruitful to examine Benders cuts that rule out neighbourhoods but outside the local-branching restrictions, i.e., non-$k$-OPT neighbourhoods. Cuts of this type could also deserve a geometric study in terms of showing that these are the deepest or that they define faces of high-dimension.

%% file: references.tex
\bibliographystyle{informs2014}

%% file: main.bbl
\begin{thebibliography}{1}
    \bibitem{ijpr2022}
    Avgerinos I., Mourtos I.,  Vatikiotis S., Zois G.,
    \newblock  Scheduling unrelated machines with job splitting, setup resources and sequence dependency
    \newblock {\em International Journal of Production Research, 61:16}, 5502-5524, (2023)
    
    \bibitem{MIM}
    Avgerinos, I., Mourtos, I., Vatikiotis, S., Zois, G.,
    \newblock Exact methods for tardiness objectives in production scheduling.
    \newblock {\em IFAC-PapersOnLine, 55:10}, 2487-2492, (2022)

    \bibitem{beck}
    Beck, C.J.,
    \newblock Checking-Up on Branch-and-Check.
    \newblock {\em CP 2010: Principles and Practice of Constraint Programming}, 84-98, (2010)

    \bibitem{asl}     Beheshti Asl N., MirHassani S.A.,     \newblock Accelerating benders decomposition: multiple cuts via multiple solutions.     \newblock {\em Journal of Combinatorial Optimization, 37}, 806–826 (2019)

    \bibitem{Bektas}
    Bekta\c s T., Hamzadayi A., Ruiz R.,
    \newblock Benders decomposition for the mixed no-idle permutation flowshop
    scheduling problem.
    \newblock {\em Journal of Scheduling , 23}, 513–523, (2020).
    
    \bibitem{Bek19}
    Bektur G., Saraç T.,
    \newblock A mathematical model and heuristic algorithms for an unrelated parallel machine scheduling problem with sequence-dependent setup times, machine eligibility restrictions and a common server.
    \newblock {\em Computers \& Operations Research, 103}, 46-63, (2019).
    
    \bibitem{benders}
    Benders J.F.,
    \newblock Partitioning procedures for solving mixed-variables programming problems. 
    \newblock {\em Numerische Mathematik, 4}, 238-252 (1962)
    
    \bibitem{Bul17}	
    Bülbül, K., Şen, H.
    \newblock An exact extended formulation for the unrelated parallel machine total weighted completion time problem.
    \newblock {\em Journal of Scheduling, 20}, 373–389, (2017).

 \bibitem{Chen09}	Chen, JF,     \newblock Scheduling on unrelated parallel machines with sequence- and machine-dependent setup times and due-date constraints.     \newblock {\em The International Journal of Advanced Manufacturing Technology, 44}, 1204–1212, (2009).
    
    
\bibitem{Dog79}	    Dogramaci A., Surkis J.     \newblock Evaluation of a Heuristic for Scheduling Independent Jobs on Parallel Identical Processors.     \newblock {\em Management Science, 25:12}, 1208-1216, (1979).
    

    \bibitem{Pey19}	
    Fanjul-Peyro, L., Ruiz, R., Perea, F.
    \newblock Reformulations and an exact algorithm for unrelated parallel machine scheduling problems with setup times.
    \newblock {\em Computers \& Operations Research, 101}, 173-182, (2019).
    
    \bibitem{Pey20}	
    Fanjul-Peyro, L.
    \newblock Models and an exact method for the Unrelated Parallel Machine scheduling problem with setups and resources.
    \newblock {\em Expert Systems with Applications: X, 5}, (2020).
    
    \bibitem{fischetti03}
    Fischetti M., Lodi A.,
    \newblock Local Branching. 
    \newblock {\em Mathematical Programming, 98}, 23-47 (2003).
    
  \bibitem{Fotakis16}
    Fotakis, D., Milis, I., Papadigenopoulos, O., Vassalos, V., Zois, G.,
    \newblock Scheduling MapReduce jobs on identical and unrelated processors.
    \newblock {\em Theory of Computing Systems, 64:5}, 754-782 (2016).	
    
    
    \bibitem{gomory}
    Gomory R.E.,
    \newblock Outline of an Algorithm for Integer Solutions to Linear Programs. 
    \newblock {\em Bulletin of the American Mathematical Society, 64}, 275-278 (1958)
    
    
    \bibitem{Hall97}	
    Hall, L.A.
    \newblock Approximation Algorithms for Scheduling.
    \newblock {\em Approximation Algorithms for NP-hard Problems, Dorit S. Hochbaum}, 1-45, (1997).
    
    \bibitem{harjunkoski}
    Harjunkoski I., Grossmann I.E.,
    \newblock Decomposition techniques for multistage scheduling problems using mixed-integer and constraint programming methods. 
    \newblock {\em Computers and Chemical Engineering, 26}, 1533-1552 (2002)
    
    \bibitem{hooker03}
    Hooker J.N., Ottoson G.,
    \newblock Logic-Based Benders Decomposition. 
    \newblock {\em Mathematical Programming, 96}, 33-60 (2003)
    
    \bibitem{hooker07}
    Hooker J.N.,
    \newblock Planning and Scheduling by Logic-Based Benders Decomposition. 
    \newblock {\em Operations Research, 55:3}, 588-602 (2007)
    
    \bibitem{hooker19}
    Hooker J.N.,
    \newblock  Logic-Based Benders Decomposition for Large-Scale Optimization.
    \newblock {\em Velásquez-Bermúdez J.M., Khakifirooz M., Fathi M., (eds) Large Scale Optimization in Supply Chains and Smart Manufacturing. Springer Optimization and Its Applications, vol 149}, (2019)
    
    
    \bibitem{jain}
    Jain V., Grossmann I.E.,
    \newblock Algorithms for hybrid MILP/CP models for a class of optimization problems. 
    \newblock {\em INFORMS Journal on Computing, 13}, 258-276 (2001)
    
    \bibitem{jeihoonian}
    Jeihoonian M., Zanjani M.K., Gendreau M.,
    \newblock Accelerating Benders decomposition for closed-loop supply chain network design: Case of used durable products with different quality levels. 
    \newblock {\em European Journal of Operational Research, 251}, 830-845 (2016)
    
    \bibitem{Lee21}
    Kim H.J., Lee J.H.,
    \newblock Scheduling uniform parallel dedicated machines with job splitting, sequence-dependent setup times, and multiple servers.
    \newblock {\em Computers \& Operations Research, 126}, (2021).
    
    
    \bibitem{Kim20}
    Kim, J.G., Song, S., Jeong, B.
    \newblock Minimising total tardiness for the identical parallel machine scheduling problem with splitting jobs and sequence-dependent setup times.
    \newblock {\em International Journal of Production Research, 58:6}, 1628-1643, (2020).
    

    \bibitem{Pinedo97}	
    Lee Y.H., Pinedo, M.,
   \newblock Scheduling jobs on parallel machines with sequence-dependent setup times.
    \newblock {\em European Journal of Operational Research, 100:3}, 464 – 474, (1997).
    
    \bibitem{Len77}	
    Lenstra, J.K., Rinnooy Kan, A.H.G., Brucker, P.
    \newblock Complexity of Machine Scheduling Problems.
    \newblock {\em Annals of Discrete Mathematics, 1}, 343-362, (1977).
    
    \bibitem{Liaw03}	
    Liaw, C.F., Lin Y.K., Cheng, C.Y., Chen, M.
    \newblock Scheduling unrelated parallel machines to minimize total weighted tardiness.
    \newblock {\em Computers \& Operations Research, 30}, 1777-1789, (2003).
    
    \bibitem{Mae20}	
    Maecker, S., Shen, L.
    \newblock Solving parallel machine problems with delivery times and ardiness objectives.
    \newblock {\em Annals of Operations Research, 285}, 315-334, (2020).
    
    \bibitem{magnati}
    Magnati T.L., Wong R.T.,
    \newblock Accelerating Benders Decomposition: Algorithmic Enhancement and Model Selection Criteria. 
    \newblock {\em Operations Research, 29:3}, 464-484 (1981)
    
    \bibitem{Pes22}	
    Pessoa, A.A., Bulhões, T., Nesello, V., Subramanian, A.
    \newblock Exact Approaches for Single Machine Total Weighted Tardiness Batch Scheduling
    \newblock {\em INFORMS Journal on Computing}, (2022).
    
    
    \bibitem{Puchinger05}
    Puchinger J., Raidl G.R.,
    \newblock  Combining Metaheuristics and Exact Algorithms in Combinatorial Optimization: A Survey and Classification.
    \newblock {\em In: Mira, J., Álvarez, J.R. (eds) Artificial Intelligence and Knowledge Engineering Applications: A Bioinspired Approach. IWINAC 2005. Lecture Notes in Computer Science, vol 3562}, 41-53 (2005)
    
    \bibitem{raidl14}
    Raidl G.R., Baumhauer T., Hu B.,
    \newblock Speeding Up Logic-Based Benders’ Decomposition by a Metaheuristic for a Bi-Level Capacitated Vehicle Routing Problem. 
    \newblock {\em In: Blesa, M.J., Blum, C., Voß, S. (eds) Hybrid Metaheuristics. HM 2014}, 183-197 (2014)
    
    \bibitem{raidl15}
    Raidl G.R., Baumhauer T., Hu B.,
    \newblock Boosting an Exact Logic-Based Benders Decomposition Approach by Variable Neighborhood Search. 
    \newblock {\em Electronic Notes in Discrete Mathematics, 47}, 149-156 (2015)
    
    \bibitem{rei09}
    Rei W., Cordeau J.F., Gendreau M., Soriano P.,
    \newblock Accelerating Benders Decomposition by Local Branching. 
    \newblock {\em INFORMS Journal on Computing, 21:2}, 333-345 (2009)
    
    
    \bibitem{seo}
    Seo K., Joung S., Lee C., Park S., 
    \newblock A Closest Benders Cut Selection Scheme for Accelerating the Benders Decomposition Algorithm. 
    \newblock {\em INFORMS Journal on Computing}, (2022)

    \bibitem{thor}
    Thorsteinsson E.S., 
    \newblock Branch-and-Check: A Hybrid Framework Integrating Mixed Integer Programming and Constraint Logic Programming. 
    \newblock {\em CP 2001: Principles and Practice of Constraint Programming}, 16-30 (2001)
    
    \bibitem{tran}
    Tran T.T., Araujo A.J., Beck C., 
    \newblock Decomposition Methods for the Parallel Machine Scheduling Problem with Setups. 
    \newblock {\em INFORMS Journal on Computing, 28:1}, 83-95 (2016)
    
    
    \bibitem{Yal00}	
    Yalaoui, F., Chu, C.,
    \newblock Parallel Machine Scheduling to Minimize Total Completion Time with Release Dates Constraints.
    \newblock {\em IFAC Proceedings Volumes, 33:17}, 715-722, (2000).
    

    \bibitem{zhu}	
    Zhu G., Bard J.F., Yu G.,
    \newblock A Branch-and-Cut Procedure for the Multimode Resource-Constrained Project-Scheduling Problem.
    \newblock {\em INFORMS Journal on Computing, 18:3}, 377-390, (2006).

\end{thebibliography}
